\documentclass{amsart}
\usepackage{amsmath,amsfonts,amsthm,amssymb,setspace}
\newtheorem{thm}{Theorem}[subsection]
\newtheorem*{mthm}{Main Theorem}
\newtheorem{lemma}{Lemma}[subsection]

\newtheorem{coro}{Corollary}[subsection]

\begin{document}

\title{The Beilinson Equivalence for Differential Operators and Lie Algebroids}

\author{Greg Muller}

\begin{abstract}
Let $\mathcal{D}$ be the ring of differential operators on a smooth irreducible affine variety $X$ over $\mathbb{C}$; or, more generally, the enveloping algebra of any locally free Lie algebroid on $X$.  The category of finitely-generated graded modules of the Rees algebra $\widetilde{\mathcal{D}}$ has a natural quotient category $qgr(\widetilde{\mathcal{D}})$ which imitates the category of modules on $Proj$ of a graded commutative ring. We show that the derived category $D^b(qgr(\widetilde{\mathcal{D}}))$ is equivalent to the derived category of finitely-generated modules of a sheaf of algebras $E$ on $X$ which is coherent over $X$.  This generalizes the usual Beilinson equivalence for projective space, and also the Beilinson equivalence for differential operators on a smooth curve used by Ben-Zvi and Nevins in \cite{BZN} to describe the moduli space of left ideals in $\mathcal{D}$.
\end{abstract}


\maketitle


\section{Introduction}

So as to appeal to a wider audience, the main theorem will be outlined in the specific case of differential operators. Then, the context will be expanded to the more general setting of Lie algebroids.

\subsection{Main Theorem: Differential Operators.}

Let $X$ be a smooth irreducible affine variety over $\mathbb{C}$ of dimension $d$, and let $\mathcal{D}$ be the ring of (algebraic) differential operators on $X$.  The ring $\mathcal{D}$ is a deformation of $Sym_X\mathcal{T}$, which is the ring of regular functions on the cotangent bundle of $X$.  Therefore, one strategy to study $\mathcal{D}$-modules is to take known methods for studying $Sym_X\mathcal{T}$-modules and see if they can be extended.

A useful technique in the study of vector bundles on a variety is to compactify them fiberwise to get a $\mathbb{P}^n$ bundle over $X$ and then use more powerful tools on $\mathbb{P}^n$ than are available in the affine case.  One such powerful tool is the Beilinson equivalence, which says that the derived category of coherent sheaves on $\mathbb{P}^n$ is equivalent to the derived category of f.g. modules of a certain quiver $Q_n$.  In the case of $\mathbb{P}^n$-bundles, there is a similar derived equivalence to an algebra which is a relative version of the quiver $Q_n$.  The purpose of this paper is to develop first the right notion of compactifying the algebra $\mathcal{D}$, and then to show there is an analog of the Beilinson equivalence.

The fiberwise compactification of $Sym_X\mathcal{T}$ is given by $Proj(\widetilde{Sym_X\mathcal{T}})$, where $\widetilde{Sym_X\mathcal{T}}$ is the Rees algebra $\oplus_{i\in \mathbb{N}}(Sym_X^{\leq i}\mathcal{T})t^i$, and where $t$ is a central variable.  The algebra $\mathcal{D}$ has a natural filtration by the degree of the operator, so we can define the Rees algebra $\widetilde{\mathcal{D}}:=\oplus_{i\in \mathbb{N}}\mathcal{D}^it^i$.  Unfortunately, there is no $Proj(\widetilde{\mathcal{D}})$, since it non-commutative.

However, there is an abelian category which imitates the category of coherent modules on $Proj(\widetilde{\mathcal{D}})$, the 'non-commutative projective geometry' of Artin and Zhang \cite{ArtinZhang}.  Let $gr(\widetilde{\mathcal{D}})$ be the category of f.g. graded $\widetilde{\mathcal{D}}$-modules, and let $tors(\widetilde{\mathcal{D}})$ denote the subcategory of graded modules non-zero in only finitely many degrees.  Then there is a quotient category $qgr(\widetilde{\mathcal{D}}):=gr(\widetilde{\mathcal{D}})/tors(\widetilde{\mathcal{D}})$, with quotient functor $\pi$.  We regard this category as the category of coherent modules on the 'non-commutative space' $Proj(\widetilde{\mathcal{D}})$.

Define the object $T:=\oplus_{i=0}^d\pi\widetilde{\mathcal{D}}(-i)\in qgr(\widetilde{\mathcal{D}})$, where $(-i)$ denotes a shift in the grading.  Then, for any $M\in D^b(qgr(\widetilde{\mathcal{D}}))$, $\mathbf{R}Hom_{qgr(\widetilde{\mathcal{D}})}(T,M)$ is a (derived) right $Hom_{qgr(\widetilde{\mathcal{D}})}(T,T)$-module.

It is not hard to show that the algebra $E:={Hom}_{qgr(\widetilde{\mathcal{D}})}(T,T)^{op}$ is given by
\[\left(\begin{array}{ccccc}
\mathcal{O}_X & \mathcal{D}^1 & \mathcal{D}^2 & \cdots & \mathcal{D}^d \\
0 & \mathcal{O}_X & \mathcal{D}^1 & \cdots & \mathcal{D}^{d-1} \\
0 & 0 & \mathcal{O}_X & \cdots & \mathcal{D}^{d-2} \\
\vdots & \vdots & \vdots & \ddots & \vdots \\
0 & 0 & 0 & \cdots & \mathcal{O}_X \\
\end{array}\right)\]
where matrix multiplication is defined by the maps $\mathcal{D}^n\otimes \mathcal{D}^m\rightarrow \mathcal{D}^{n+m}$.  Therefore, $\mathbf{R}{Hom}_{qgr(\widetilde{\mathcal{D}})}(T,-)$ defines a functor from $D^b(qgr(\widetilde{\mathcal{D}}))$ to $D^b(mod(E))$, where $mod(E)$ is the category of sheaves of f.g. left $E$-modules.

\begin{mthm} (The Beilinson Equivalence for Differential Operators)
The functor $\mathbf{R}{Hom}_{qgr(\widetilde{\mathcal{D}})}(T,-):D^b(qgr(\widetilde{\mathcal{D}}))\rightarrow D^b(mod(E))$ is an equivalence of triangulated categories.
\end{mthm}

The proof is in two parts.  The first is to show that every object in $D^b(qgr(\widetilde{\mathcal{D}}))$ can be resolved by summands of $T$; this will follow from developing an $\mathcal{O}_X$-relative version of Koszul duality theory.  The second is to show that ${Ext}^i_{qgr(\widetilde{\mathcal{D}})}(T,T)=0$ for all $i>0$.  This will establish $T$ as a compact generator of $qgr(\widetilde{\mathcal{U}})$ with derived endomorphism algebra $E$, and the theorem will follow from the usual tilting theorems in this case.

The chief advantage of this result is that it relates the representation theory of $\mathcal{D}$ to the representation theory of $E$, an algebra which is finitely generated over $\mathcal{O}_X$, and is much easier to study.  This idea was used by Ben-Zvi and Nevins \cite{BZN} when $X$ is a curve to relate isomorphism classes of left ideals in $\mathcal{D}$ to specific collections of data in $D^b(mod(X))$.  This was used to characterize the moduli space of such ideals, generalizing earlier results of Berest and Wilson \cite{YuriWilsonIdeals}.

\subsection{The Main Theorem: Lie Algebroids.}

Lie algebroids are a simultaneous generalization of rings of differential operators and of Lie algebras.  Studying them can be very useful for understanding those aspects of the representation theory of Lie algebras which have an analogous statement for the representation theory of differential operators.  However, there are many interesting Lie algebroids which are neither Lie algebras nor differential operators.

A Lie algebroid $L$ on $X$ is a coherent $\mathcal{O}_X$-module which is also a sheaf of Lie algebras, together with a map (the \textbf{anchor map}) $\tau:L\rightarrow\mathcal{T}_X$ which is a map of $\mathcal{O}_X$-modules and of sheaves of Lie algebras.\footnote{Here and throughout the paper, $\mathcal{T}_X$ or $\mathcal{T}$ will denote the tangent bundle of $X$ over $\mathbb{C}$, which is the same as the module of $\mathbb{C}$-linear derivations on $\mathcal{O}_X$.}  However, instead of the bracket on $L$ commuting with the $\mathcal{O}_X$ multiplication, it fails to commute in a way controlled by the anchor map: $[l,fl']=f[l,l']+d_{\tau(l)}(f)\cdot l'$, where $l,l'\in L$, $f\in \mathcal{O}_X$ and $d_{\tau(l)}$ is the derivative along $\tau(l)$.  Note that $\mathcal{T}_X$ itself is a Lie algebroid with the identity map as anchor, and it is via this Lie algebroid that this section generalizes the previous section.  From now on, we will require that $L$ is locally-free as a $\mathcal{O}_X$-module (which $\mathcal{T}_X$ is), and its rank will be denoted by $n$.

$L$ has a notion of a representation on a $\mathcal{O}_X$-module, and there is a corresponding universal enveloping algebra $\mathcal{U}_XL$.  For $L=\mathcal{T}_X$, the enveloping algebra is the ring of differential operators $\mathcal{D}$.  $\mathcal{U}_XL$ has a natural filtration, and the Rees algebra $\widetilde{\mathcal{U}_XL}$ can be defined as the graded algebras $\oplus_{i\in \mathbb{N}}\mathcal{U}^i_XL\cdot t^i$, where $t$ is a central variable. The categories $gr(\widetilde{\mathcal{U}_XL})$, $tors(\widetilde{\mathcal{U}_XL})$, and $qgr(\widetilde{\mathcal{U}_XL})$ all have identical definitions to the previous case.

Let $T:=\oplus_{i=0}^{n}\pi(\mathcal{U}_XL(-i))$.  The algebra $E:={Hom}_{qgr(\widetilde{\mathcal{U}_XL})}(T,T)^{op}$ is
\[\left(\begin{array}{ccccc}
\mathcal{O}_X & \mathcal{U}_X^1L & \mathcal{U}^2_XL & \cdots & \mathcal{U}^n_XL \\
0 & \mathcal{O}_X & \mathcal{U}^1_XL & \cdots & \mathcal{U}^{n-1}_XL \\
0 & 0 & \mathcal{O}_X & \cdots & \mathcal{U}^{n-2}_XL \\
\vdots & \vdots & \vdots &\ddots & \vdots \\
0 & 0 & 0 & \cdots & \mathcal{O}_X \\
\end{array}\right)\]
Then the functor $\mathbf{R}{Hom}_{qgr(\widetilde{\mathcal{U}_XL})}(T,-)$ goes from $D^b(qgr(\widetilde{U}_XL))$ to $D^b(mod(E))$.

\begin{mthm} (The Beilinson Equivalence for Lie Algebroids)
The functor $\mathbf{R}{Hom}_{qgr(\widetilde{\mathcal{U}_XL})}(T,-):D^b(qgr(\widetilde{\mathcal{U}_XL}))\rightarrow D^b(mod(E))$ is an equivalence of triangulated categories.
\end{mthm}

It is this theorem that will actually be proved; it will specialize to the first main theorem in the case of $L=\mathcal{T}_X$.  However, this generality allows other interesting corollaries, some of which were already well known.  The Beilinson equivalence for $\mathbb{P}^n$ and for $\mathbb{P}^n$-bundles is reproved (but not in a substantially different way), as well as a Beilinson equivalence for the quantum spaces of Lie algebras introduced by Le Bruyn and van den Bergh in \cite{LeBruynvdB}.  Section \ref{section:examples} addresses these examples, as well as some applications of this theorem.

The techniques of this paper can also be readily applied to a wider array of examples than this paper covers (Section \ref{subsection:nonexamples} describes what properties of $\widetilde{\mathcal{U}}$ are necessary to make the main theorem work).  It should also be noted that the results of this paper hold in the case of non-affine $X$, provided the Beilinson functor is defined correctly; the necessary proofs for this are contained in Appendix \ref{section:local}.

\subsection{Outline of the Paper.}

Section 2 contains (without proof) the basics of Lie algebroids and non-commutative projective geometry which will be used throughout the paper.  Section 3 is a rather technical tangent which builds up the notion of tensor product in non-commutative projective geometry, far enough to define Fourier-Mukai transforms.  Section 4 develops the Koszul theory of $\widetilde{\mathcal{U}}$ and produces several kinds of canonical resolutions.  Section 5 contains the proof of the main theorem.  Section 6 concludes the body of the paper by outlining some basic examples of interest, some quick applications, and explores how the scope of the theorem can be generalized.  Appendix A collects several important computations regarding the quadratic dual algebra $\widetilde{\mathcal{U}}^\perp$, which are necessary for certain proofs.  Appendix B deals with the case of non-affine $X$, as well as proving the naturality of the constructions in this paper with respect to localization.

\subsection{Acknowledgements.}

This paper is most immediately indebted to \cite{BZN} for inspiration, as well as to David Ben-Zvi and Thomas Nevins for many accomodating emails explaining the more delicate details of their paper.  Another paper which had a significant formative affect on this paper is \cite{KKO}, which this paper borrows techniques from quite liberally.

On a personal level, the author is especially thankful to Yuri Berest, whose tutelage and support have been invaluable throughout this paper's preparation.  Many other people were gracious in sharing both their conversation and insights; an almost certainly incomplete list includes David Ben-Zvi, Thomas Nevins, Peter Samuelson, Frank Moore, Margaret Bruns, and Tomoo Matsumura.

\section{Preliminaries.}


\subsection{Notation for Graded Modules.}  First, let us fix some notation.  Let $A$ be a graded algebra, and let $M$ be a graded left $A$-module.  Let $M(i)$ denote the $n$th shifting functor, so that $[M(i)]_j=M_{i+j}$.  If $N$ is a graded left $A$-module, then $Hom_{Gr}(M,N)$ will denote the degree zero maps from $M$ to $N$, while $\underline{Hom}_A(M,N)$ will denote $\oplus_{i\in \mathbb{Z}} Hom(M,N(i))$.  If $N$ is a graded right $A$-module, then $N{\otimes}_AM$ will denote the graded tensor, which is a graded vector space, while $N\circledcirc_A M$ will denote $(N\otimes_AM)_0$, the degree zero part.

\subsection{Non-Commutative Projective Geometry.}

If $A$ is a non-commutative algebra, then there is no general consensus as to what sort of object $Spec(A)$ should be, or even if can exist at all.  However, instead of trying to build a locally ringed space to call $Spec(A)$, we can simply work with the category $Mod(A)$, thought of as the category of quasi-coherent sheaves on the non-existent $Spec(A)$.  Since most questions one might ask about a scheme can be restated as a question about its category of modules, this allows many questions of a geometric flavor to be answered.

Now, if $A$ is a positively-graded algebra, we can similarly bypass the need for a space $Proj(A)$ and instead concern ourselves with its category of modules. We use the projective Serre equivalence as a recipe for what this category should be.  Let $Gr(A)$ be the category of graded left $A$-modules, and let $Tors(A)$ be the full subcategory of modules such that, for every $m\in T\in Tors(A)$, $A_{\geq n}\cdot m=0$ for some $n$.  Let $gr(A)$ denote the category of finitely generated graded left modules, and $tors(A):=gr(A)\cap Tors(A)$.  If $A$ is commutative, then the category of quasi-coherent sheaves on $Proj(A)$ is equivalent to the quotient category $QGr(A):=Gr(A)/Tors(A)$, while the category of coherent modules is equivalent to $qgr(A):=gr(A)/tors(A)$.

For not-necessarily commutative $A$, we will think of $QGr(A):=Gr(A)/Tors(A)$ as the category of quasi-coherent modules on the undefined space $Proj(A)$.  This perspective was first put forward by Artin and Zhang in \cite{ArtinZhang}.  We collect the necessary facts about $QGr(A)$ below without proofs, with their page listings in \cite{ArtinZhang}.
\begin{itemize}
\item The shifting functor descends to a functor on $QGr(A)$, and so $\underline{Hom}_{QGr(A)}$ is well defined.
\item (pg. 235) $QGr(A)$ has enough injectives.
\item (pg. 234) The quotient functor $\pi:Gr(A)\rightarrow QGr(A)$ is exact.
\item (pg. 234) The functor $\pi$ has a right adjoint $\omega:QGr(A)\rightarrow Gr(A)$ which is left exact.  Because 
    \[\omega(M)=\underline{Hom}_{Gr(A)}(A, \omega(M))=Hom_{Qgr(A)}(\pi A,M)\] $\omega(M)$
    should be regarded as the 'graded global section functor'.  In this vein, $\mathbf{R}^i\omega(M)$ is the analog of the $i$th graded cohomology of $M$.
\item (pg. 234) $\pi\omega(M)=M$.
\item (pg. 234) If $A$ is left noetherian, then $\omega\pi(M)=\lim_{\rightarrow} \underline{Hom}_{Gr(A)}(A_{\geq n},M)$.
\item (pg. 233) Every module $M\in Gr(A)$ has a maximal submodule $\tau(M)$ in $Tors(A)$.  It can be explicitly defined by
    \[\tau(M):=\lim_{\rightarrow}\underline{Hom}_{Gr(A)}(A/A_{\geq n},M)\]
    where the left $A$-module structure on $\tau(M)$ comes from the right $A$-module structure on $A/A_{\geq n}$. It is a left exact functor, and its derived functors $\mathbf{R}^i\tau(M)$ coincide with the $i$th local cohomology of $M$ at the ideal $A_{\geq 1}$, at least when $A$ is generated in degree 0 and 1.
\item (pg. 241) The defining inclusion $\tau(M)\hookrightarrow M$ and the adjunction map $M\rightarrow \omega\pi(M)$ fit together to give an exact triangle in $D(Gr(A))$:
    \[\mathbf{R}\tau(M)\rightarrow M \rightarrow \mathbf{R}\omega\pi(M)\rightarrow \mathbf{R}\tau(M)[1]\]
    In particular, we have an exact sequence
    \[0\rightarrow \tau(M)\rightarrow M\rightarrow \omega\pi(M)\rightarrow \mathbf{R}^1\tau(M)\rightarrow 0\]
    and natural equivalences $\mathbf{R}^i\omega\pi(M)\simeq \mathbf{R}^{i+1}\tau(M)$ for $i\geq 1$.
\item (pg. 243) A module $M\in Gr(A)$ is said to satisfy the \textbf{$\chi_i$-condition} if, for all $d$ and all $j\leq i$, there is an $n_0$ such that for all $n\geq n_0$, $\underline{Ext}^i_{Gr(A)}(A/A_{\geq n},M)_{\geq j}$ is a finitely-generated $A$ module.  $M$ has the \textbf{$\chi$-condition} if it satisfies $\chi_i$ for all $i$.
\item (pg. 273) (Serre Finiteness) Let $A$ be left noetherian and satisfy $\chi$, and let $M\in gr(A)$.  Then, for all $i\geq 1$, $\mathbf{R}^i\omega\pi(M)_d$ is a finitely generated $A_0$-module for all $d$, and is zero if $d$ is sufficiently large.
\end{itemize}

In the case of $A$ left noetherian, there is a more useful definition of $\mathbf{R}\omega\pi(M)$.
\begin{lemma}\label{Lemma:Cech Lemma}
Let $A$ be left noetherian.  For $M\in Gr(A)$, there is an isomorphism in $D(Gr(A))$:
\[\mathbf{R}\omega\pi(M)\simeq \mathbf{R}\omega\pi(A){\otimes}_A^\mathbf{L} M\]
\end{lemma}
\begin{proof}
This follows immediately from the isomorphisms $\mathbf{R}\underline{Hom}_{Gr(A)}(A_{\geq n},M)\simeq \mathbf{R}\underline{Hom}_{Gr(A)}(A_{\geq n},A){\otimes}_A^\mathbf{L} M$.
\end{proof}
Applying this for $M=\mathbf{R}\omega\pi(A)$,
\begin{coro}\label{Coro:Cech Idempotent}
There is an isomorphism in the derived category:
\[\mathbf{R}\omega\pi(A){\otimes}^\mathbf{L}_A\mathbf{R}\omega\pi(A)\simeq \mathbf{R}\omega\pi(\mathbf{R}\omega\pi(A))=\mathbf{R}\omega\pi(A)\]
\end{coro}

\subsection{Lie Algebroid Basics.}

The study of Lie algebroids comes from the infinitesmal study of Lie groupoids.\footnote{Hence the name.  It has nothing to do with objects that are more properly called algebroids (at least not when these were named); it is a pun on 'groupoid'.}  However, much like Lie algebras, Lie algebroids are intrinsically interesting, even without a corresponding Lie groupoid in mind.  For a more detailed reference, consult \cite{Mackenzie}.

An (algebraic) \textbf{Lie algebroid} on $X$ is an $\mathcal{O}_X$-module $L$ with
\begin{itemize}
\item a \textbf{Lie bracket} on $L$ which makes it into a Lie algebra.
\item an \textbf{anchor map}, an $\mathcal{O}_X$-module map $\tau:L\rightarrow \mathcal{T}$.
\end{itemize}
The bracket and the $\mathcal{O}_X$-module structure on $L$ are \emph{not} necessarily compatible in the simplest way\footnote{The simplest way would be that each are defined arrow-theoretically in the category of the other; this would be an $\mathcal{O}_X$-Lie algebra.  They correspond to Lie algebroids with trivial anchor map.}; instead, the bracket and the $\mathcal{O}_X$-multiplication satisfy the relation:
\[ [l,al']=a[l,l']+d_{\tau(l)}(a)\cdot l'\]
One consequence of this relation is that $\mathcal{O}_X\oplus L$ becomes a Lie algebras by the bracket $[(r,l),(r',l')]=(d_{\tau(l)}(r')-d_{\tau(l')}(r),[l,l'])$.  In this paper, the only Lie algebroids which will be considered are those such that $L$ is a locally-free coherent $\mathcal{O}_X$-module; this will be assumed from here on.

A Lie algebroid comes with instructions on how to commute two sections of $L$ past each other (the bracket) and how to commute sections of $L$ past sections of $\mathcal{O}_X$ (the anchor).  This naturally leads to the consideration of the universal algebra generated by $L$ and $\mathcal{O}_X$ which obey the given commutation relations.  Let $\mathcal{U}_XL$ be the quotient of the universal enveloping algebra of the Lie algebra $\mathcal{O}_X\oplus L$ by the relations $(1,0)=1$ and $(a,0)\otimes (a',l)=(aa',al)$ ($1$ the unit, $a\in\mathcal{O}_X$, and $l\in L$); this is called the \textbf{universal enveloping algebra} of $L$.  The algebra $\mathcal{U}_XL$ is the central object of study in this paper; for simplicity, it will be denoted $\mathcal{U}$ when $X$ and $L$ are clear.

$\mathcal{O}_X$ has a canonical structure of a left $\mathcal{U}$-module, by the action $a\cdot a'=aa'$ and $l\cdot a=d_{\tau(l)}(a)$ for $a,a'\in \mathcal{O}_X$ and $l\in L$.  The 'action on 1' map $\mathcal{U}\rightarrow \mathcal{O}_X$ which sends $\sigma$ to $\sigma\cdot 1$ is a left $\mathcal{U}$-module map which presents $\mathcal{O}_X$ as a quotient of $\mathcal{U}$ as a left module over itself.  Note however, that there is no canonical right $\mathcal{U}$-module structure on $\mathcal{O}_X$.

$\mathcal{U}$ is naturally filtered by letting the image of $\mathcal{O}_X$ be degree 0 and the image of $L$ be degree 1.  The subspace $\mathcal{U}^1$ is a (not-necessarily central) $\mathcal{O}_X$-bimodule which fits into a short exact sequence of $\mathcal{O}_X$-bimodules:
\[0\rightarrow \mathcal{O}_X\rightarrow \mathcal{U}^1\rightarrow L\rightarrow 0\]
The Rees algebra $\widetilde{\mathcal{U}}$ is the graded algebra defined as $\oplus_{i\in\mathbb{N}} \mathcal{U}^i\cdot t^i$, where $t$ is a central element.  The Rees algebra has the property that $\widetilde{\mathcal{U}}/(t-1)=\mathcal{U}$.  It can also be defined directly as a quotient of the tensor algebra $T_X \mathcal{U}^1$ by the relation $\partial\otimes \partial' -\partial'\otimes \partial =[\partial,\partial']\otimes t$, where $\partial,\partial'\in\mathcal{U}^1$ and $t$ denotes $1\in \mathcal{O}_X\subset\mathcal{U}^1$ (as opposed to the unit of the algebra).

$\widetilde{\mathcal{U}}/t$ is equal to $\overline{\mathcal{U}}$, the \textbf{associated graded algebra} of $\mathcal{U}$ which is usually defined as $\oplus_{i\in \mathbb{N}}\mathcal{U}^i/\mathcal{U}^{i-1}$.  Because the commutator of a degree $i$ element and a degree $j$ element in $\mathcal{U}$ is of degree at most $i+j-1$, the algebra $\overline{\mathcal{U}}$ is commutative.  In fact, $\overline{\mathcal{U}}$ is isomorphic to $Sym_XL$, the symmetric algebra generated by $L$ (this is the PBW theorem for Lie algebroids).  This is also isomorphic to $f^*(\mathcal{O}_{L^*})$, the total space of the dual bundle to $L$ pushed forward along the bundle map $f:L^*\rightarrow X$.

A nice consequence of the PBW theorem for Lie algebroids is that $\widetilde{\mathcal{U}}^i$ is projective and finitely-generated as both a left and right $\mathcal{O}_X$-module (though not as a bimodule). This is because $\widetilde{\mathcal{U}}^i/\widetilde{\mathcal{U}}^{i-1}=\overline{\mathcal{U}}^i$ is f.g. projective, and so $\widetilde{\mathcal{U}}^i$ has a finite composition sequence consisting entirely of f.g. projectives.  As a consequence, $\widetilde{\mathcal{U}}$ satisfies the $\chi$-condition, and so $\mathbf{R}^i\omega\pi(M)_j$ is a finitely generated $\mathcal{O}_X$-module for all $i$ and $j$, as long as $M$ is finitely generated as a $\widetilde{\mathcal{U}}$-module.

\section{Tensoring and Fourier-Mukai Transforms.}

We need to generalize an important technique from commutative projective geometry to the non-commutative setting; that of the Fourier-Mukai transform.  Let $X$ be a scheme, and let $K$ be any module on $X\times X$, or more generally any derived object in $D^b(Mod(X\times X))$ (equivalently, $K$ is a derived $\mathcal{O}_X$-bimodule).  Given any $M\in D^b(Mod(X))$, $K\otimes_X^{\mathbf{L}}M\in D^b(Mod(X\times X))$, and so it can be pushed forward along the projection $p_1:X\times X\rightarrow X$ onto the first factor to give $\mathbf{R}p_{1*}(K\otimes_X^{\mathbf{L}}M)\in D^b(Mod(X))$.  The functor $M\rightarrow \mathbf{R}p_{1*}(K\otimes_X^{\mathbf{L}} M)$ is called the \textbf{Fourier-Mukai transform of $K$}.  These have been studied extensively, for references check \cite{HuybrechtsFourierMukai}.
\subsection{Tensor Products.}

For $A$ a positively-graded algebra, the categories $Gr(A)$ and $gr(A)$ don't have a tensor product in the sense of a bifunctorial map $Gr(A)\times Gr(A)\rightarrow Gr(A)$.  The tensor product here is a bifunctorial map ${\otimes}_A:Gr(A^{op})\times Gr(A)\rightarrow Gr(\mathbb{C})$.  Subsequently taking the degree zero part gives a map ${\circledcirc}_A:Gr(A^{op})\times Gr(A)\rightarrow Vect$.

Naively, one would hope that this descends to some kind of map ${\circledcirc}_A:QGr(A^{op})\times QGr(A)\rightarrow Vect$.  However, for this to descend to a map on quotient categories, we would need that $T{\circledcirc}_A M=M'{\circledcirc}_A T'=0$ for $T\in Tors(A^{op})$ and $T'\in Tors(A)$.  This is just not true; take, for example, $A_0{\otimes}_A A$ or $A{\otimes}_A A_0$, which are both isomorphic to $A_0$ as a vector space.

So, instead of trying to push the multiplication forward along $\pi$, we can pull the multiplication back along $\omega$.  Given $\pi{M}\in QGr(A^{op})$ and $\pi{N}\in QGr(A)$, define \[\pi{M}{\circledcirc}_{A}\pi{N}:= \omega \pi{M}{\circledcirc}_A\omega\pi{N}=(\omega \pi{M}{\otimes}_A\omega\pi{N})_0\]
As a combination of left and right exact functors, this will not in general be a left or right exact bifunctor.  Nonetheless, the natural related derived construction is $(\mathbf{R}\omega\pi{M}{\otimes}^\mathbf{L}_A\mathbf{R}\omega\pi{N})_0$ (for $\pi{M}\in D^b(QGr(A^{op}))$ and $\pi{N}\in D^b(QGr(A))$).

\subsection{The Category of Quotient Bimodules.} The point of these tensoring constructions is to be able to define the Fourier-Mukai transforms on this category; however, we still need to know where the kernels of the transforms live.  Let $A^e:= A\otimes A^{op}$, which has the property that $A^{op}$-modules are the same as $A$-bimodules; note that it is a bigraded algebra.  Let $\mathbb{G}r(A^e)$ be the category of bigraded $A^e$-modules, which is the same as the category of bigraded $A$-bimodules.  Let $\mathbb{T}ors(A^e)$ be the subcategory of $\mathbb{G}r(A^e)$ such that, for every $m\in T\in \mathbb{T}ors(A^e)$, there is some $n$ such that $A_{\geq n}mA_{\geq n}=0$.  Let $Q\mathbb{G}r(A^e)$ denote the quotient category $\mathbb{G}r(A^e)/\mathbb{T}ors(A^e)$.

$Q\mathbb{G}r(A^e)$ satisfies all the same properties that were listed for $QGr(A)$, or at least analogous properties.\footnote{The key properties that make this work are that $\mathbb{T}ors(A^e)$ is a dense subcategory, every object in $\mathbb{G}r(A^e)$ has a maximal torsion-submodule, and that $\mathbb{G}r(A^e)$ has enough injectives.  See [AZ] pg. 234, or [Po] Sect. 4.4.}  The only difference is the structure of the functors $\omega$ and $\tau$, which may be given by (where $\underline{Hom}$ now denotes a bigraded $Hom$)
\[\omega\pi(M):=\lim_{n\rightarrow}\underline{Hom}_{\mathbb{G}r(A^e)}(A_{\geq n}\otimes A_{\geq n},M)\]
\[\tau(M):= \lim_{n\rightarrow} \underline{Hom}_{\mathbb{G}r(A^e)}((A\otimes A) /(A_{\geq n}\otimes A_{\geq n}),M)\]

In certain nice cases, the derived functor $\mathbf{R}\omega\pi$ has a more useful definition.
\begin{lemma}\label{Lemma:Split Cech Lemma}
Let $A$ be left and right noetherian.  For $M\in \mathbb{G}r(A^e)$, there is an isomorphism in $D(\mathbb{G}r(A^e))$:
\[\mathbf{R}\omega\pi(M)\simeq \mathbf{R}\omega\pi(A){\otimes}_A^\mathbf{L} M{\otimes}_A^\mathbf{L} \mathbf{R}\omega\pi(A)\]
\end{lemma}
\begin{proof}
Consider the directed system $A_{\geq m}\otimes A_{\geq m'}$, as $m$ and $m'$ run over the integers, with the maps being the natural inclusions.  This directed system has a sub-directed system $A_{\geq n}\otimes A_{\geq n}$ which is coinitial, in the sense that any object $A_{\geq m}\otimes A_{\geq m'}$ has a surjection from some $A_{\geq n}\otimes A_{\geq n}$ (for instance, $n=\min(m,m')$).  Therefore, there is an isomorphism of direct limits:
\[\lim_{n\rightarrow} \mathbf{R}\underline{Hom}_{\mathbb{G}r(A^e)}(A_{\geq n}\otimes A_{\geq n},M)\simeq\lim_{m\rightarrow}\lim_{m'\rightarrow} \mathbf{R}\underline{Hom}_{\mathbb{G}r(A^e)}(A_{\geq m}\otimes A_{\geq m'},M)\]
By adjunction, this second $\mathbf{R}\underline{Hom}$ becomes \[\lim_{m\rightarrow}\lim_{m'\rightarrow}\mathbf{R}\underline{Hom}_{Gr(A)}(A_{\geq m},\mathbf{R}\underline{Hom}_{Gr(A^{op})}(A_{\geq m'},M))\]
\[=\lim_{m\rightarrow}\lim_{m'\rightarrow}\mathbf{R}\underline{Hom}_{Gr(A)}(A_{\geq m},A){\otimes}^\mathbf{L}_A \mathbf{R}\underline{Hom}_{Gr(A^{op})}(A_{\geq m'},M)\]
\[=\lim_{m\rightarrow}\lim_{m'\rightarrow}\mathbf{R}\underline{Hom}_{Gr(A)}(A_{\geq m},A){\otimes}^\mathbf{L}_A M{\otimes}^\mathbf{L}_A\mathbf{R}\underline{Hom}_{Gr(A^{op})}(A_{\geq m'},A)\]
The last two equalities use that $A_{\geq m}$ is noetherian as a left and right $A$-module.
This final expression is then equal to $\mathbf{R}\omega\pi(A){\otimes}^\mathbf{L}_A M{\otimes}^\mathbf{L}_A \mathbf{R}\omega\pi(A)$.
\end{proof}

\subsection{Fourier-Mukai Transforms.} Now, given any object $K\in D^b(\mathbb{G}r(A^e))$, define the derived functor $F_{K}$ on $D^b(QGr(A))$ by:
\[F_{K}(\pi{M}):=  \pi(\mathbf{R}\omega\pi({K}){\otimes}^\mathbf{L}_A \mathbf{R}\omega(\pi{M}))_{\bullet,0}\]
This has a simpler form for nice $A$.

\begin{lemma}
If $A$ is left and right noetherian, then
\[F_{K}(\pi M)=\pi(K\circledcirc_A^\mathbf{L} \mathbf{R}\omega\pi(M))=\pi(\mathbf{R}\omega\pi(K)\circledcirc_A^\mathbf{L} M)\]
\end{lemma}
\begin{proof}
By Lemma \ref{Lemma:Split Cech Lemma} and Lemma \ref{Lemma:Cech Lemma}, this is equal to
\[\pi(\mathbf{R}\omega\pi(A){\otimes}_A^\mathbf{L} {K} {\otimes}_A^\mathbf{L}\mathbf{R}\omega\pi(A) {\otimes}^\mathbf{L}_A \mathbf{R}\omega\pi(A){\otimes}_A^\mathbf{L}M)_{\bullet,0}\]
By Corollary \ref{Coro:Cech Idempotent}, this is
\begin{equation}\label{Eqn:FMeq}\pi(\mathbf{R}\omega\pi(A){\otimes}_A^\mathbf{L} {K} {\otimes}_A^\mathbf{L} \mathbf{R}\omega\pi(A){\otimes}_A^\mathbf{L}M)_{\bullet,0}\end{equation}
Applying Lemma \ref{Lemma:Cech Lemma} twice and using that $\pi\mathbf{R}\omega\pi=\pi$, this is equal to
\[\pi(\mathbf{R}\omega\pi(A){\otimes}_A^\mathbf{L} {K} {\otimes}_A^\mathbf{L} \mathbf{R}\omega\pi({M}))_{\bullet,0}=\pi(\mathbf{R}\omega\pi(K\otimes_A^\mathbf{L}\mathbf{R}\omega\pi({M})))
=\pi(K\circledcirc_A^\mathbf{L}\mathbf{R}\omega\pi({M}))\]
Instead, we could apply Lemma \ref{Lemma:Split Cech Lemma} to equation \ref{Eqn:FMeq} to get
\[\pi(\mathbf{R}\omega\pi(K){\otimes}_A^\mathbf{L} M)_{\bullet,0}=\pi(\mathbf{R}\omega\pi(K)\circledcirc_A^\mathbf{L}M)\]
This concludes the proof.
\end{proof}

Given any exact triangle $A\rightarrow B\rightarrow C\rightarrow A[1]$ in $D^b(\mathbb{G}r(A^e)$, there is an associated exact triangle of functors $F_A\rightarrow F_B\rightarrow F_C\rightarrow F_A[1]$, in the sense that for any $\pi{M}\in D^b(QGr(A^e))$, there is an exact triangle:
\[ F_A(\pi{M})\rightarrow F_B(\pi{M})\rightarrow F_C(\pi{M})\rightarrow F_A(\pi{M})[1]\]
Therefore, a functor $F_K$ may be resolved by other, simpler functors by resolving $\pi K$ into simpler objects in $D^b(\mathbb{G}r(A^e))$.

\subsection{The Identity Functor as a Fourier-Mukai Transform.}

Even the identity functor on $D^b(QGr(A))$ arises as a Fourier-Mukai transform.  Let $\widetilde{\Delta}$ be the bigraded $A$-bimodule such that $\widetilde{\Delta}_{i,j}=A^{i+j}$, where $A^{k}=0$ in negative degrees.  $\widetilde{\Delta}$ has the property that $\widetilde{\Delta}\circledcirc_A M=(\widetilde{\Delta}{\otimes}_AM)_{\bullet,0}=M$ for all $M\in Gr(A)$.  As an immediate corollary, $\widetilde{\Delta}$ is flat as a right $A$-module.  If $A$ is noetherian, the Fourier-Mukai transform $F_{\widetilde{\Delta}}(\mathcal{M})$ is $\pi(\widetilde{\Delta}\circledcirc_A \mathbf{R}\omega(\mathcal{M}))$, which is $\pi(\mathbf{R}\omega(\mathcal{M}))=\mathcal{M}$.  Therefore, $F_{\widetilde{\Delta}}$ is the identity functor.

However, $\widetilde{\Delta}$ is not the only object in $\mathbb{G}r(A^e)$ whose corresponding Fourier-Mukai transform is the identity.  After all, all that matters is the image under $\pi$ in $Q\mathbb{G}r(A^e)$.  Let $\Delta$ be the bigraded $A$-bimodule such that $\Delta_{i,j}=A^{i+j}$ when $i\geq0$ and $j\geq 0$, and zero otherwise.  There is a natural inclusion $\Delta\hookrightarrow \widetilde{\Delta}$, and $(\widetilde{\Delta}/\Delta)_{i,j}=0$ if $i\geq 0$ and $j\geq 0$.  $\omega\pi(\widetilde{\Delta}/\Delta)=0$, and so $\pi(\Delta)=\pi(\widetilde{\Delta})$.  Then, the Fourier Mukai transform $F_{\Delta}$ is also the identity.

The point of this is now that producing a resolution of $\Delta$ in $\mathbb{G}r(A^e)$ will give a resolution of the identity, which in turn will give a resolution of any object.

\section{Koszul Duality for Lie Algebroids.}\label{section:koszul}

The goal now is to show that any object in $QGr(\widetilde{\mathcal{U}})$ can be resolved by sums of objects of the form $\pi\widetilde{\mathcal{U}}(-i)$, for $0\leq i\leq d$.  This will be accomplished by developing the Koszul theory for the algebra $\mathcal{U}$ over $X$.  The two main results of this will be:
\begin{itemize}
\item A canonical projective resolution of $\mathcal{O}_X$ as a left $\widetilde{\mathcal{U}}$-module, called the \textbf{left Koszul resolution}.
\item For any $\pi M\in QGr(\widetilde{\mathcal{U}})$, a projective resolution of $\pi M$, called the \textbf{Beilinson resolution}.
\end{itemize}
%
%
%

The key observation is that the definition of the universal enveloping algebra gives a surjective map $T_X\mathcal{U}^1\rightarrow \widetilde{\mathcal{U}}$, whose kernel is generated by elements of degree 2 in $T_X\mathcal{U}$.  A \textbf{relatively quadratic algebra over $X$} is an algebra with a surjective map from $T_XB$ for some $\mathcal{O}_X$-bimodule $B$, whose kernel is generated in degree 2.  The motivating case is that of quadratic algebras over $Spec(\mathbb{C})$ (the non-relative case), where there is an elaborate theory of Koszul resolutions and duality.  This section is an extension of those techniques to the current setting.

\subsection{The Quadratic Dual Algebra.}

Let $R$ be the $\mathcal{O}_X$-bimodule which is the kernel of the map $\mathcal{U}^1\otimes_X\mathcal{U}^1\rightarrow \mathcal{U}^2$.  Note that $R$ is the degree 2 part of the kernel of $T_X\mathcal{U}^1\rightarrow \widetilde{\mathcal{U}}$, which generates the whole kernel as a two-sided ideal.  By the definition of the universal enveloping algebra, this is the $\mathcal{O}_X$-bimodule generated by $\partial\otimes \partial'-\partial'\otimes \partial -[\partial,\partial']\otimes 1$ for $\partial,\partial'\in \mathcal{U}^1$.  

From now on, for $M$ any right $\mathcal{O}_X$-module, let $M^*$ denote the left $\mathcal{O}_X$-module $Hom_{-X}(M,\mathcal{O}_X)$ (as right $\mathcal{O}_X$-modules)\footnote{$Hom_{-X}$ will denote the $Hom$ as right $\mathcal{O}_X$-modules, when there is also a left $\mathcal{O}_X$-structure.  Similarly, $Hom_{X-}$ will denote the $Hom$ as left $\mathcal{O}_X$-modules.}; analogously, for $M$ any left $\mathcal{O}_X$-module, let $^*M$ denote the right $\mathcal{O}_X$-module $Hom_{X-}(M,\mathcal{O}_X)$.  When $M$ is a $\mathcal{O}_X$-bimodule, $M^*$ and $^*M$ are also $\mathcal{O}_X$-bimodules, which are potentially non-isomorphic.

Let $J^i$ be $\,^*(\mathcal{U}^i)$ , which is called the \textbf{bimodule of $i$-jets}.  Since the $\mathcal{U}^i$ are finitely generated and projective as right $\mathcal{O}_X$-modules, there is an isomorphism $\,^*(\mathcal{U}^1\otimes_X\mathcal{U}^1)\simeq \,^*(\mathcal{U}^1)\otimes_X \,^*(\mathcal{U}^1)\simeq J^1\otimes_X J^1$.  The map $\mathcal{U}^1\otimes_X\mathcal{U}^1\rightarrow \mathcal{U}^2$ then induces an inclusion $J^2\hookrightarrow J^1\otimes_X J^1$, which can be characterized as the subset of right $\mathcal{O}_X$-module maps $\mathcal{U}^1\otimes_X \mathcal{U}^1\rightarrow \mathcal{O}_X$ which kill $R\subset \mathcal{U}^1\otimes_X \mathcal{U}^1$.

Now, let $\widetilde{\mathcal{U}}^\perp$ denote the quotient of the tensor algebra $T_XJ^1$ by the two-sided ideal generated by $J^2$ as sitting inside the degree 2 part.  The algebra $\widetilde{\mathcal{U}}^\perp$ is called the \textbf{Koszul dual}, or the \textbf{quadratic dual} algebra.\footnote{Note that we have made an asymmetric choice, in looking at the dual of $\mathcal{U}^1$ as a left $\mathcal{O}_X$-module, rather than as a right $\mathcal{O}_X$-module.  Then, perhaps, this should be called the left Koszul dual.  This choice was motivated by the fact that $J^1$ has much nicer properties than $(\mathcal{U}^1)^*$, which results in a nicer presentation of $\widetilde{\mathcal{U}}^\perp$. However, the right Koszul dual algebra would still have been sufficient for the purposes of this paper.}  In Appendix A.2, it is shown that $\widetilde{\mathcal{U}}^\perp=\underline{Ext}^{\bullet}_{\widetilde{\mathcal{U}}-}(\mathcal{O}_X,\mathcal{O}_X)$, where $\underline{Ext}$ is the graded $Ext$; arguably, this is where it draws all its important properties.

Appendix A contains proofs of several interesting results about $\widetilde{\mathcal{U}}^{\perp}$; however, we only need the following facts.
\begin{itemize}
\item (Corollary \ref{coro:dual fg proj}) For all $i$, $\widetilde{\mathcal{U}}^{\perp i}$ is finitely-generated projective as a left and right $\mathcal{O}_X$-module.
\item (Corollary \ref{coro:dual finite}) For $i>n+1$ or $i<0$, $\widetilde{\mathcal{U}}^{\perp i}=0$.
\item (Corollary \ref{coro:dual duality}) Let $\omega_L$ be $\Lambda^n_XL^*$, the top exterior power of
the dual to $L$.  Then $(\widetilde{\mathcal{U}}^{\perp i})^*= \omega_L^*\otimes_X \widetilde{\mathcal{U}}^{\perp (n+1-i)}$ and $\,^*(\widetilde{\mathcal{U}}^{\perp i})= \widetilde{\mathcal{U}}^{\perp (n+1-i)}\otimes_X\omega_L^* $.
\end{itemize}


\subsection{The Left Koszul Complex.}
The Koszul dual algebra now lets us construct a canonical resolution of the left $\widetilde{\mathcal{U}}$-module $\mathcal{O}_X$, which will be important for the coming steps.

The multiplication map $m_{\widetilde{\mathcal{U}}^\perp}: \widetilde{\mathcal{U}}^{\perp i-1}\otimes_X J^1\rightarrow \widetilde{\mathcal{U}}^{\perp i}$ induces a right dual map
\[m^\vee_{\mathcal{U}^\perp}:(\widetilde{\mathcal{U}}^{\perp i})^*\rightarrow (\widetilde{\mathcal{U}}^{\perp i-1}\otimes_X J^1)^*\simeq (J^1)^*\otimes_X(\widetilde{\mathcal{U}}^{\perp i-1})^*\simeq \mathcal{U}^1\otimes_X(\widetilde{\mathcal{U}}^{\perp i-1})^*\]
Define a composition map,
\[k^i:\widetilde{\mathcal{U}}(-i)\otimes_X (\widetilde{\mathcal{U}}^{\perp i})^*\rightarrow \widetilde{\mathcal{U}}(-i)\otimes_X \mathcal{U}^1\otimes_X (\widetilde{\mathcal{U}}^{\perp i-1})^*\rightarrow \widetilde{\mathcal{U}}(-i+1)\otimes_X (\widetilde{\mathcal{U}}^{\perp i-1})^*\]
where the first map is the above map $m_{\mathcal{U}^\perp}^\vee$, and the second map is the multiplication map $m_{\mathcal{U}}:\widetilde{\mathcal{U}}(-i)\otimes_X \mathcal{U}^1\rightarrow \widetilde{\mathcal{U}}(-i+1)$.  Let $K^i_{(X,L)}$ (or $K^i$ when $X$ and $L$ are clear) denote the left $\widetilde{\mathcal{U}}$-module $\widetilde{\mathcal{U}}(-i)\otimes_X (\widetilde{\mathcal{U}}^{\perp i})^*$. Note that $K^i=0$ if $i<0$ or $i>n+1$.
\begin{thm}
The map $k^i:K^i\rightarrow K^{i-1}$ makes $K^\bullet$ into a complex of left $\widetilde{\mathcal{U}}$-modules called the \textbf{left Koszul complex}.
\end{thm}
\begin{proof}
The square of the Koszul boundary, $k^2$, is $m_{\mathcal{U}}m^\vee_{\mathcal{U}^\perp}m_{\mathcal{U}}m^\vee_{\mathcal{U}^\perp}$.  However, the middle two maps can be commuted, since they involve disjoint terms in the tensor product.  Therefore, $k^2=(m_{\mathcal{U}})^2(m_{\mathcal{U}^\perp}^\vee)^2$, which is the composition
\[\widetilde{\mathcal{U}}(-i)\otimes_X (\widetilde{\mathcal{U}}^{\perp i})^*\rightarrow \widetilde{\mathcal{U}}(-i)\otimes_X \mathcal{U}^1\otimes_X \mathcal{U}^1\otimes_X (\widetilde{\mathcal{U}}^{\perp i-2})^* \rightarrow \widetilde{\mathcal{U}}(-i+2)\otimes_X (\widetilde{\mathcal{U}}^{\perp i-2})^*\]
The map $(m_{\mathcal{U}^\perp}^\vee)^2$ is the map
\[Hom_{X-}(\widetilde{\mathcal{U}}^{\perp i},\mathcal{O}_X)\rightarrow Hom_{X-}(\widetilde{\mathcal{U}}^{\perp i-2}\otimes_X J^1\otimes_X J^1,\mathcal{O}_X)\]
right dual to multiplication.  Everything in the image of this map necessarily kills $\widetilde{\mathcal{U}}^{\perp i-2}\otimes_X J^2\subset \widetilde{\mathcal{U}}^{\perp i-2}\otimes_X J^1\otimes_X J^1$, which translates to the image of $(m_{\mathcal{U}^\perp}^\vee)^2$ being contained in $R\otimes_X (\widetilde{\mathcal{U}}^{\perp i-2})^*$.  Then, it is clear that the multiplication map $(m_{\mathcal{U}})^2$ kills anything in $\widetilde{\mathcal{U}}(-i)\otimes_X R\otimes_X (\widetilde{\mathcal{U}}^{\perp i-2})^*$.  Therefore, $k^2=0$.
\end{proof}
The construction of the left Koszul complex commutes with localization in the natural way, as per the following lemma.
\begin{lemma}
Let $X'\subset X$ be an open subscheme of $X$, and $L'$ the localization of $L$.  Then the left Koszul complex $K^\bullet_{(X',L')}$ of the Lie algebroid $(X',L')$ is equal to the localization of the left Koszul complex $K^\bullet_{(X,L)}$.
\end{lemma}\begin{proof}
This is Lemma \ref{lemma:koszul local} in the Appendix.
\end{proof}
We are finally ready for the most meaningful fact about the Koszul complex, that it resolves $\mathcal{O}_X$ as a left $\widetilde{\mathcal{U}}$-module.
\begin{thm}
The natural quotient map $K^0=\widetilde{\mathcal{U}}\rightarrow \mathcal{O}_X$ makes $K^\bullet$ into a resolution of $\mathcal{O}_X$; that is, the complex $K^\bullet$ is exact in positive degrees, and its cohomology in degree zero is exactly the image of the augmentation map.
\end{thm}
\begin{proof}
The strategy of the proof will be a succession of cases of increasing generality.

\begin{itemize}
\item \textbf{$X=Spec(\mathbf{k})$ ($\mathbf{k}$ a field), $L$ abelian.} This is the classical case of Koszul duality for $Sym_\mathbf{k}L$ and $\Lambda_{\mathbf{k}}L$.  A proof can be found in \cite{Weibel}, page 114.
\item \textbf{$X$ a regular local ring, $L$ abelian.} Because $X$ is local, $L$ being projective implies that it is free, specifically that $L=\mathcal{O}_X\otimes_\mathbf{k}L/m$ where $\mathbf{k}$ is the residue field.  The Rees algebra $\widetilde{\mathcal{U}}$ is isomorphic to the symmetric algebra $Sym_X L=\mathcal{O}_X\otimes_\mathbf{k}Sym_\mathbf{k}L/m$.  The quadratic dual algebra $\widetilde{\mathcal{U}}^\perp$ is then the corresponding exterior algebra $Alt_X L^*=\mathcal{O}_X\otimes_\mathbf{k} Alt_\mathbf{k} L^*/m$.  The left Koszul complex $K^\bullet_{(X,L)}$ is then $\mathcal{O}_X\otimes_\mathbf{k} K^\bullet_{(Spec(\mathbf{k}),L/m)}$, where $L/m$ is the Lie algebroid restricted to the residue field $\mathbf{k}$.  Since the theorem is true for $K^\bullet_{(Spec(\mathbf{k}),L/m)}$ by the previous case, it is then true here.
\item \textbf{$X$ arbitrary, $L$ abelian.} Let $\pi:X_p\rightarrow X$ be the localization at some prime $p$, and let $L_p=\pi^* L$.  By the lemma before the theorem, $\pi^*K^\bullet_{(X,L)}=K^\bullet_{(X_p,L_p)}$.  Since localization is exact, we have that \[\pi^*H^i\left(K^\bullet_{(X,L)}\right)=H^i\left(\pi^*K^\bullet_{(X,L)}\right)\]
    The two facts together imply that $\pi^*H^i\left(K^\bullet_{(X,L)}\right)=H^i\left(K^\bullet_{(X_p,L_p)}\right)$.  The previous case of the theorem implies that this second group vanishes for $i>0$, and is isomorphic to $\mathcal{O}_{X_p}$ for $i=0$. Since this fact is true at any prime $p$, it is true everywhere, and so the theorem is true.
\item \textbf{$X$ arbitrary, $L$ arbitrary.} Consider a family of Lie algebroids $(X,L_\hbar)$, $\hbar\in \mathbb{C}$, where the bracket $[-,-]_\hbar:= \hbar[-,-]$ and $\tau_\hbar:=\hbar\tau$.  In this notation, $L_1$ is the original Lie algebroid, and $L_0$ is the Lie algebroid with zero bracket and anchor.  Notice that, for $\hbar\neq 0$, the Lie algebroid $L_\hbar$ is isomorphic to $L=L_1$, by the scaling map $l\rightarrow \hbar^{-1}l$ for $l\in L$.  Therefore, the Koszul complex $K^\bullet_{(X,L_\hbar)}$ with parameter $\hbar$ is isomorphic to $K^\bullet_{(X,L)}$ away from $\hbar=0$.  By the previous case of the theorem, the theorem is true for $L_0$.  Since the cohomology must be constant in some neighborhood of $\hbar=0$ and the complexes are isomorphic for all other $\hbar$, the theorem is true for all $\hbar$, in particular $\hbar=1$.
\end{itemize}
\end{proof}

Since $\mathcal{U}^{\perp i}$ is a f.g. projective right $\mathcal{O}_X$-module, then $(\mathcal{U}^{\perp i})^*$ is a f.g. projective left $\mathcal{O}_X$-module.  Therefore, $K^i$ is a projective left $\widetilde{\mathcal{U}}$-module, and the left Koszul resolution is a projective resolution of $\mathcal{O}_X$ as a $\widetilde{\mathcal{U}}$-module.

%
%

There is also a right Koszul complex $K^\bullet_{right}$ whose terms are $(\widetilde{\mathcal{U}}^{\perp i} )^*\otimes_X \widetilde{\mathcal{U}}(-i)$, with boundary right dual to the multiplication map $\mathcal{U}^1\otimes_X \widetilde{\mathcal{U}}^{i-1}\rightarrow \widetilde{\mathcal{U}}^i$.  This is again a projective resolution of $\mathcal{O}_X$, this time as a right $\widetilde{\mathcal{U}}$-module.  The proofs are analogous.

\subsection{The Koszul Bicomplex.}

The next step is to combine the left and right Koszul complexes into a Koszul bicomplex, which can then be used to extract a resolution of the diagonal.

Let $\mathbb{K}^{i,j}$ be the $\widetilde{\mathcal{U}}$-bimodule $\widetilde{\mathcal{U}}(-i)\otimes_X (\widetilde{\mathcal{U}}^{\perp (i+j)})^*\otimes_X \widetilde{\mathcal{U}}(-j)$.  The left Koszul boundary map acts on the first two terms, and sends $\mathbb{K}^{i,j}$ to $\mathbb{K}^{i-1,j}$; the right Koszul boundary map acts on the last two terms, and sends $\mathbb{K}^{i,j}$ to $\mathbb{K}^{i,j-1}$.
\begin{lemma}
These two boundary maps, $k_{left}$ and $k_{right}$, make $\mathbb{K}^{i,j}$ into a bicomplex of $\widetilde{\mathcal{U}}$-bimodules called the \textbf{Koszul bicomplex} (making sure to obey the Koszul sign rule for commuting odd-degree maps).
\end{lemma}
\begin{proof}
It is immediate that the two boundaries square to zero themselves.  Thus, all that remains to check is that $(k_{left}+k_{right})$ squares to zero, which by the Koszul sign rule is equivalent to $k_{left}$ and $k_{right}$ commuting.

Since multiplication in $\widetilde{\mathcal{U}}^{\perp}$ is associative, the multiplication map $J^1\otimes_X \widetilde{\mathcal{U}}^{\perp i-2}\otimes_X J^1\rightarrow \widetilde{\mathcal{U}}^{\perp i}$ doesn't depend on the order of multiplication.  Dualizing gives the desired fact that $k_{left}$ and $k_{right}$ commute.
\end{proof}

The terms of the Koszul bicomplex are bigraded $\widetilde{\mathcal{U}}$-bimodules, and so an element in this complex can have a graded bidegree (it's bigrading as a $\widetilde{\mathcal{U}}$-bimodule) as well as a homological bidegree (which term of the bicomplex it is in).  The space of elements with graded bidegree $(p,q)$ and homological bidegree $(i,j)$ will be denoted $\mathbb{K}^{i,j}_{p,q}$, and it is equal to $\mathcal{U}^{p-i}\otimes_X (\widetilde{\mathcal{U}}^{\perp (i+j)})^*\otimes_X \mathcal{U}^{q-j}$.


\subsection{The Resolution of the Diagonal.}

We can now produce a resolution of the diagonal. Define the complex $\mathcal{K}_\Delta$ to be such that $\mathcal{K}_\Delta^i=\ker (d_r: \mathbb{K}^{i,0}\rightarrow \mathbb{K}^{i,-1})$, together with the boundary $d_l$ inherited from $\mathbb{K}$.  Because $\mathbb{K}^{0,-1}=0$, we have that $\mathcal{K}_{\Delta}^0=\mathbb{K}^{0,0}=\widetilde{\mathcal{U}}\otimes_X \widetilde{\mathcal{U}}$.

Recall from the previous section the diagonal object $\Delta\in \mathbb{G}r(\widetilde{\mathcal{U}}^e)$, a bigraded $\widetilde{\mathcal{U}}$-bimodule.  There is a canonical surjection $\widetilde{\mathcal{U}}\otimes_X \widetilde{\mathcal{U}}\rightarrow \Delta$, which in bidegree $(p,q)$ is the multiplication map $\mathcal{U}^p\otimes_X\mathcal{U}^q\rightarrow \mathcal{U}^{p+q}$.

\begin{thm}
The canonical surjection $\mathcal{K}_{\Delta}\rightarrow \Delta$ makes $\mathcal{K}_{\Delta}$ into a resolution of $\Delta$.  Accordingly, the complex $\mathcal{K}_{\Delta}$ is called a \textbf{resolution of the diagonal}.
\end{thm}
\begin{proof}
First, we show that the map $\mathcal{K}^0_\Delta\rightarrow \Delta$ gives an augmentation of the complex; that is, it kills the image of $\mathcal{K}^1_\Delta$ in $\mathcal{K}^0_\Delta$.  By definition, $\mathcal{K}^1_\Delta$ is the kernel of \[\widetilde{\mathcal{U}}(-1)\otimes_X \mathcal{U}^1\otimes_X \widetilde{\mathcal{U}}\rightarrow\widetilde{\mathcal{U}}(-1)\otimes_X \widetilde{\mathcal{U}}(1)\] This map is given by multiplying the last two terms.  However, since the composition map $\mathcal{K}^1_\Delta\rightarrow \mathcal{K}^0_\Delta\rightarrow \Delta$ is given by multiplying all the terms of $\mathcal{K}^1_\Delta$ together, and because multiplication in $\mathcal{U}$ is associative, this composition must be zero.

Now, define the \textbf{truncated Koszul bicomplex} $\widehat{\mathbb{K}}^{i,j}$ to be equal to $\mathbb{K}^{i,j}$ when $j\geq 0$, and $0$ otherwise.  For a fixed graded bidegree $(p,q)$, the term $\widehat{\mathbb{K}}^{i,j}_{p,q}$ vanishes for $i>p$, $j>q$ or $i+j<0$.  Therefore, in any fixed graded bidegree, the bicomplex $\widehat{\mathbb{K}}$ is bounded.  This means that both the horizontal-then-vertical spectral sequence and the vertical-then-horizontal spectral sequence converge to total cohomology of $\widehat{\mathbb{K}}$.

Taking horizontal cohomology first, the rows are all right Koszul complexes tensored with $\widetilde{\mathcal{U}}$, and so we get
\[E_1^{i,j}=\left\{\begin{array}{cc}
\widetilde{\mathcal{U}}(j)\otimes_X \mathcal{O}_X(-j) & if\; j=-i\geq 0 \\
0 & otherwise \\
\end{array}\right\}\]
Therefore, the spectral sequence collapses on the first page, and we have \[H^0(Tot(\widehat{\mathbb{K}}))=\bigoplus_{j=0}^{\infty}\widetilde{\mathcal{U}}(j)\otimes_X \mathcal{O}_X(-j),\; H^{\neq 0}(Tot(\widehat{\mathbb{K}}))=0\]

Taking vertical cohomology first, the rows are either left Koszul complexes tensored with $\widetilde{\mathcal{U}}$, or they are left Koszul complexes which have been brutally truncated.  Therefore,
\[E_1^{i,j}=\left\{\begin{array}{cc}
 \mathcal{O}_X(j)\otimes_X\widetilde{\mathcal{U}}(-j) & if\; j=-i\geq 1 \\
\mathcal{K}^{i}_\Delta & if \; i\geq 0, j=0\\
0 & otherwise \\
\end{array}\right\}\]
Therefore, the spectral sequence collapses on the second page, and we have
\[ H^0(Tot(\widehat{\mathbb{K}})) = H^0(\mathcal{K}_\Delta)\oplus \left(\bigoplus_{j=1}^\infty \mathcal{O}_X(j)\otimes_X \widetilde{\mathcal{U}}(-j)\right),\; H^{i\neq 0}(Tot(\widehat{\mathbb{K}}))=H^i(\mathcal{K}_{\Delta})\]

Comparing the two results, $\mathcal{K}_{\Delta}$ is exact outside degree zero, and we have that
\[H^0(\mathcal{K}_{\Delta})\oplus \left(\bigoplus_{j=1}^\infty \mathcal{O}_X(j)\otimes_X \widetilde{\mathcal{U}}(-j)\right)
=\bigoplus_{j=0}^{\infty}\widetilde{\mathcal{U}}(j)\otimes_X \mathcal{O}_X(-j)\]
Looking in graded bidegree $(p,q)$, we have that $H^0(\mathcal{K}_\Delta)=\mathcal{U}^{p+q}$ if and only if $p,q\geq 0$.  Therefore, the map $H^0(\mathcal{K}_\Delta)\rightarrow \Delta$ induced by the augmentation is an isomorphism.
\end{proof}

The power of this theorem comes from the structure of $\mathcal{K}_{\Delta}$.  To see this structure, define $\Omega^i_R$ to be the kernel of the $i$-th boundary in the right Koszul complex: \[d_r:(\widetilde{\mathcal{U}}^{\perp i})^*\otimes_X \widetilde{\mathcal{U}}(-i)\rightarrow (\widetilde{\mathcal{U}}^{\perp i-1})^*\otimes_X \widetilde{\mathcal{U}}(-i+1)\]
Since $\widetilde{\mathcal{U}}^{\perp j}=0$ for $j>n+1$, $\Omega^j_R=0$ for $j>n$.
It is clear from the definition of $\mathcal{K}_{\Delta}$ that $ \mathcal{K}_{\Delta}^i=\widetilde{\mathcal{U}}(-i)\otimes_X \Omega^i_R(i)$.
\begin{coro} The resolution of the diagonal then has the form:
\[\Delta\leftarrow \widetilde{\mathcal{U}}\otimes_X \widetilde{\mathcal{U}}\leftarrow \widetilde{\mathcal{U}}(-1)\otimes_X \Omega^1_R(1) \leftarrow ... \leftarrow \widetilde{\mathcal{U}}(-i)\otimes_X \Omega^i_R(i) \leftarrow ...\leftarrow \widetilde{\mathcal{U}}(-n)\otimes_X \Omega^n_R(n)\]
\end{coro}

There is a mirror image version of this, where $\mathcal{K}_{\Delta}$ is replaced by $ker(d_l:\mathbb{K}^{0,i}\rightarrow \mathbb{K}^{-1,i})$.  Defining \[\Omega_L^i:=ker\left(d_l:\widetilde{\mathcal{U}}(-i)\otimes_X (\widetilde{\mathcal{U}}^{\perp i})^*\rightarrow
\widetilde{\mathcal{U}}(-i+1)\otimes_X (\widetilde{\mathcal{U}}^{\perp i-1})^*\right),\]
all the same arguments work to show that the following is also a resolution of the diagonal:
\[\Delta\leftarrow \widetilde{\mathcal{U}}\otimes_X \widetilde{\mathcal{U}}\leftarrow \Omega^1_L(1) \otimes_X \widetilde{\mathcal{U}}(-1) \leftarrow ... \leftarrow \Omega^i_L(i)\otimes_X \widetilde{\mathcal{U}}(-i) \leftarrow ...\leftarrow  \Omega^n_L(n)\otimes_X\widetilde{\mathcal{U}}(-n)\]

\subsection{The Beilinson Resolution.}

The resolution of the diagonal then gives a resolution for every object $\pi M$ in $QGr(\widetilde{\mathcal{U}})$.

\begin{thm}
Every object $\pi(M)\in QGr(\widetilde{\mathcal{U}})$ has a resolution of the form:
\[\pi\left(\widetilde{\mathcal{U}}\otimes_X^\mathbf{L} \mathbf{R}\omega\pi(M)_0\right)\leftarrow ... \pi\left(\widetilde{\mathcal{U}}(-i)\otimes_X^\mathbf{L} \left(\Omega^i_R(i)\circledcirc_{\widetilde{\mathcal{U}}}^\mathbf{L} \mathbf{R}\omega\pi(M)\right)\right)\leftarrow ...\]
\end{thm}
\begin{proof}
The resolution of the diagonal gives a complex of Fourier-Mukai transforms.  Applying each of these to some $\pi{M}\in QGr(\widetilde{\mathcal{U}})$, we get
\[F_{\Delta}(\pi{M})\leftarrow F_{\widetilde{\mathcal{U}}\otimes_X \widetilde{\mathcal{U}}}(\pi{M})\leftarrow ...\leftarrow F_{\widetilde{\mathcal{U}}(-i)\otimes_X \Omega^i_R(i)}(\pi{M})\leftarrow ... F_{\widetilde{\mathcal{U}}(-n)\otimes_X \Omega^n_R(n)}(\pi{M})\]
The first object is $\pi{M}$, by the design of $\Delta$.  The Fourier-Mukai transform is \[F_{\widetilde{\mathcal{U}}(-i)\otimes_X \Omega^i_R(i)}(\pi M)=\pi(\mathbf{R}\omega\pi(\widetilde{\mathcal{U}}(-i)\otimes_X \Omega^i_R(i)){\otimes}_{\widetilde{\mathcal{U}}}^\mathbf{L}\mathbf{R}\omega\pi(M))_{\bullet,0}\]
By Lemma \ref{Lemma:Split Cech Lemma},
\[=\pi\left( \mathbf{R}\omega\pi(\widetilde{\mathcal{U}}){\otimes}_{\widetilde{\mathcal{U}}}^\mathbf{L}
\left(\widetilde{\mathcal{U}}(-i)\otimes_X^L\Omega^i_R(i)\right)
{\otimes}_{\widetilde{\mathcal{U}}}^\mathbf{L}\mathbf{R}\omega\pi(\widetilde{\mathcal{U}})
{\otimes}_{\widetilde{\mathcal{U}}}^\mathbf{L} \mathbf{R}\omega\pi(M)\right)_{\bullet,0}\]
which simplifies to
\[\pi\left( \mathbf{R}\omega\pi(\widetilde{\mathcal{U}}(-i))\otimes_X^\mathbf{L}\left( \Omega^i_R(i) \circledcirc_{\widetilde{\mathcal{U}}}^\mathbf{L} \mathbf{R}\omega\pi(M)\right)\right)=\pi\left(\widetilde{\mathcal{U}}(-i)\otimes_X^\mathbf{L}\left( \Omega^i_R(i) \circledcirc_{\widetilde{\mathcal{U}}}^\mathbf{L} \mathbf{R}\omega\pi(M)\right)\right)\]
\end{proof}

Note that $\left(\Omega^i_R(i)\circledcirc_{\widetilde{\mathcal{U}}}^\mathbf{L}\mathbf{R}\omega\pi(M)\right)$ is a derived object in $D^b(Coh(X))$.  Since $X$ is affine, every object in $Coh(X)$ can be resolved by copies of the structure sheaf $\mathcal{O}_X$.  This means that $\widetilde{\mathcal{U}}(-i)\otimes_X^\mathbf{L} \left(\Omega^i_R(i)\circledcirc_{\widetilde{\mathcal{U}}}^\mathbf{L}\mathbf{R}\omega\pi(M)\right)$ is quasi-isomorphic to a complex consisting of copies of $\widetilde{\mathcal{U}}(-i)$.  Therefore,
\begin{coro}\label{coro:gen T}
Every object $\pi M\in QGr(\widetilde{\mathcal{U}})$ has a resolution consisting of sums of the objects $\pi\widetilde{\mathcal{U}}$, $\pi\widetilde{\mathcal{U}}(-1)$, ... $\pi\widetilde{\mathcal{U}}(-n)$.
\end{coro}

\section{The Beilinson Equivalence.}

The previous section proved that any $\pi M\in qgr(\widetilde{\mathcal{U}})$ has a finite resolution by finite sums of the objects $\pi \widetilde{\mathcal{U}}$, $\pi\widetilde{\mathcal{U}}(-1)$, ... and $\pi\widetilde{\mathcal{U}}(-n)$. Therefore, there is always a surjection $T^{\oplus i}\rightarrow \pi M$ for large enough $i$; and so $T$ is called a \textbf{generator} for the category $qgr(\widetilde{\mathcal{U}})$.  The next question is the structure of $\mathbf{R}{Hom}_{qgr(\widetilde{\mathcal{U}})}(T,T)$.

\subsection{The Relative Gorenstein Property.}

The vanishing of the higher $Ext$'s from $T$ to itself will follow from the following property, which should be regarded as relative version of the Gorenstein property for graded algebras.  Recall that $\omega_L:=\Lambda^n_XL^*$.

\begin{lemma} (The relative Gorenstein property)
\[\underline{Ext}^i_{\widetilde{\mathcal{U}}-}(\mathcal{O}_X,\widetilde{\mathcal{U}})
=\left\{\begin{array}{cc} \omega_L(n+1) & i=n+1 \\
0 & otherwise \\
\end{array} \right.\]
\end{lemma}
\begin{proof}
Resolve $\mathcal{O}_X$ by the left Koszul resolution $K^\bullet$.  Using Corollary \ref{coro:dual duality}, which says that $(\widetilde{\mathcal{U}}^{\perp i})^*= \omega_L^*\otimes_X\widetilde{\mathcal{U}}^{\perp (n+1-i)}$ and $\,^*(\widetilde{\mathcal{U}}^{\perp i})= \widetilde{\mathcal{U}}^{\perp (n+1-i)}\otimes_X \omega_L^*$,
\begin{eqnarray*} \underline{Hom}_{\widetilde{\mathcal{U}}-}(K^i,\widetilde{\mathcal{U}})
&=& \underline{Hom}_{\widetilde{\mathcal{U}}-}(\widetilde{\mathcal{U}}(-i)\otimes_X (\widetilde{\mathcal{U}}^{\perp i})^*,\widetilde{\mathcal{U}})\\
&=& \underline{Hom}_{\widetilde{\mathcal{U}}-}(\widetilde{\mathcal{U}}(-i)\otimes_X \omega_L^*\otimes_X \widetilde{\mathcal{U}}^{\perp (n+1-i)},\widetilde{\mathcal{U}})\\
&=& Hom_{X-}(\omega_L^*\otimes_X\widetilde{\mathcal{U}}^{\perp (n+1-i)},\mathcal{O}_X)\otimes_X \widetilde{\mathcal{U}}(i)\\
&=& \,^*(\widetilde{\mathcal{U}}^{\perp (n+1-i)})\otimes_X\omega_L\otimes_X  \widetilde{\mathcal{U}}(i)\\
&=& \omega_L\otimes_X (\widetilde{\mathcal{U}}^{\perp (n+1-i)})^*\otimes_X \widetilde{\mathcal{U}}(i)
\end{eqnarray*}
Since the duality map is adjoint to the multiplication map, the boundary map on this complex is the right Koszul differential.  Therefore,
\[\mathbf{R}\underline{Hom}_{\widetilde{\mathcal{U}}-}(\mathcal{O}_X,\widetilde{\mathcal{U}})= \omega_L(n+1)[-n-1]\otimes_X K_{right}^\bullet\]
Since $K^\bullet_{right}$ is a resolution of $\mathcal{O}_X$, the theorem follows.
\end{proof}

\subsection{The Derived Endomorphism Algebra of $T$.}

The relative Gorenstein property is the key lemma in computing the structure of $\mathbf{R}{Hom}_{qgr(\widetilde{\mathcal{U}})}(T,T)$.  We then have the following lemma.
\begin{thm}\label{thm:end T}
For $i>0$, $Ext^i_{qgr(\widetilde{\mathcal{U}})}(T,T)=0$, and
\[Hom_{qgr(\widetilde{\mathcal{U}})}(T,T)
=\left(\begin{array}{ccccc}
\mathcal{O}_X & \mathcal{U}^1 & \mathcal{U}^2 & \cdots & \mathcal{U}^n \\
0 & \mathcal{O}_X & \mathcal{U}^1 & \cdots & \mathcal{U}^{n-1} \\
0 & 0 & \mathcal{O}_X & \cdots & \mathcal{U}^{n-2} \\
\vdots & \vdots & \vdots &\ddots & \vdots \\
0 & 0 & 0 & \cdots & \mathcal{O}_X \\
\end{array}\right)\]
\end{thm}
\begin{proof}
Replacing $T=\oplus_{i=0}^n\pi\widetilde{\mathcal{U}}(-i)$ gives that \begin{eqnarray*}
\mathbf{R}Hom_{qgr(\widetilde{\mathcal{U}})}(T,T)
&=&\mathbf{R}Hom_{qgr(\widetilde{\mathcal{U}})}(\oplus_{i=0}^n\pi\widetilde{\mathcal{U}}(-i),\oplus_{i=0}^n\pi\widetilde{\mathcal{U}}(-i))\\
&=& \bigoplus_{0\leq i,j \leq n}\mathbf{R}Hom_{qgr(\widetilde{\mathcal{U}})}(\pi\widetilde{\mathcal{U}}(-i),\pi\widetilde{\mathcal{U}}(-j))\\
&=& \bigoplus_{0\leq i,j \leq n} \mathbf{R}Hom_{qgr(\widetilde{\mathcal{U}})}(\pi\widetilde{\mathcal{U}},\pi\widetilde{\mathcal{U}}(i-j))\\
&=& \bigoplus_{0\leq i,j \leq n} \mathbf{R}Hom_{gr(\widetilde{\mathcal{U}})}(\widetilde{\mathcal{U}},\mathbf{R}\omega\pi\widetilde{\mathcal{U}}(i-j))\\
&=& \bigoplus_{0\leq i,j \leq n}
[\mathbf{R}\omega\pi(\widetilde{\mathcal{U}})]_{j-i}
\end{eqnarray*}
The derived object $\mathbf{R}\omega\pi(\widetilde{\mathcal{U}})$ fits into an exact triangle in $D^b(gr(\widetilde{\mathcal{U}}))$
\[\mathbf{R}\tau(\widetilde{\mathcal{U}}) \rightarrow \widetilde{\mathcal{U}} \rightarrow \mathbf{R}\omega\pi(\widetilde{\mathcal{U}})\rightarrow \]
However, the relative Gorenstein property can be used to show that $\mathbf{R}\tau(\widetilde{\mathcal{U}})$ vanishes above graded degree $-n-1$, and outside cohomological degree $n+1$.
\begin{lemma}\label{lemma: tau vanishing}
$(\mathbf{R}^k\tau(\widetilde{\mathcal{U}}))_j=0$ if $j> -n-1$ or if $k\neq n+1$.
\end{lemma}
\begin{proof}
For any $i$, there is a short exact sequence of $\widetilde{\mathcal{U}}$-modules:
\[0\rightarrow \mathcal{U}^i(-i)\rightarrow \widetilde{\mathcal{U}}/\widetilde{\mathcal{U}}^{\geq i+1}\rightarrow \widetilde{\mathcal{U}}/\widetilde{\mathcal{U}}^{\geq i}\rightarrow 0\]
where $\mathcal{U}^i(-i)$ is the left $\mathcal{O}_X$-module $\mathcal{U}^i$ concentrated in degree $i$, and  given a $\widetilde{\mathcal{U}}$-module structure by allowing $\widetilde{\mathcal{U}}^{\geq 1}$ to act trivially.  Applying $\underline{Hom}_{gr(\widetilde{\mathcal{U}})}(-,\widetilde{\mathcal{U}})$ to this sequence gives an exact triangle of derived objects
\[ \mathbf{R}\underline{Hom}_{\widetilde{\mathcal{U}}}(\widetilde{\mathcal{U}}/\widetilde{\mathcal{U}}^{\geq i},\widetilde{\mathcal{U}})\rightarrow \mathbf{R}\underline{Hom}_{\widetilde{\mathcal{U}}}(\widetilde{\mathcal{U}}/\widetilde{\mathcal{U}}^{\geq i+1},\widetilde{\mathcal{U}})\rightarrow \mathbf{R}\underline{Hom}_{\widetilde{\mathcal{U}}}(\mathcal{U}^i(-i),\widetilde{\mathcal{U}})\rightarrow\]

Note that $\mathbf{R}\underline{Hom}_{\widetilde{\mathcal{U}}}(\mathcal{U}^i(-i),\widetilde{\mathcal{U}})
=\,^*(\mathcal{U}^i)\otimes_X\mathbf{R}\underline{Hom}_{\widetilde{\mathcal{U}}}(\mathcal{O}_X,\widetilde{\mathcal{U}})(i)$.
From the relative Gorenstein property, $(\mathbf{R}^k\underline{Hom}_{\widetilde{\mathcal{U}}}(\mathcal{O}_X,\widetilde{\mathcal{U}}))_j=0$  for $j>-n-1$ or for $k\neq n+1$, and in these cases, the above triangle implies an isomorphism for all $i$ \[(\mathbf{R}^k\underline{Hom}_{\widetilde{\mathcal{U}}}(\widetilde{\mathcal{U}}/\widetilde{\mathcal{U}}^{\geq i},\widetilde{\mathcal{U}}))_j\simeq \mathbf{R}^k\underline{Hom}_{\widetilde{\mathcal{U}}}(\widetilde{\mathcal{U}}/\widetilde{\mathcal{U}}^{\geq i+1},\widetilde{\mathcal{U}})_j\]
However, in the case of $i<0$, these $\mathbf{R}^k\underline{Hom}$'s vanish, and so they vanish for all $i$.  Therefore, \[(\mathbf{R}^k\tau(\widetilde{\mathcal{U}}))_j=\lim_{i\rightarrow \infty}\mathbf{R}^k\underline{Hom}_{\widetilde{\mathcal{U}}}(\widetilde{\mathcal{U}}/\widetilde{\mathcal{U}}^{\geq i},\widetilde{\mathcal{U}})_j=0\]
\end{proof}

It immediately follows that $\mathbf{R}\omega\pi(\widetilde{\mathcal{U}})_k\simeq \widetilde{\mathcal{U}}^k=\mathcal{U}^k$ for $k\geq -n$, and so
\[\mathbf{R}Hom_{qgr(\widetilde{\mathcal{U}})}(T,T)=\bigoplus_{0\leq i,j\leq n}\mathcal{U}^{j-i}\]
Therefore, the higher $Ext$s vanish completely, and the endomorphism algebra of $T$ is given by the above algebra.
\end{proof}

\subsection{Equivalence to $D^b(Mod(E))$.}

Let $E$ denote $(Hom_{qgr}(T,T))^{op}$, the opposite algebra.  Now, given any $\pi M\in qgr(\widetilde{\mathcal{U}})$, $\mathbf{R}Hom_{qgr}(T,\pi M)$ has a right action by $Hom_{qgr}(T,T)$ by composition, and so it is a left $E$ module.  In this way, the functor $\mathbf{R}Hom_{qgr}(T,-)$ defines a functor from $D^b(qgr(\widetilde{\mathcal{U}}))$ to $D^b(mod(E))$.

This functor can be expressed in terms of the functor $\mathbf{R}\omega\pi$.  After all, as derived right $\mathcal{O}_X$-modules,
\begin{eqnarray*}
\mathbf{R}Hom_{qgr(\widetilde{\mathcal{U}})}(T,\pi M)
&=&\mathbf{R}Hom_{qgr(\widetilde{\mathcal{U}})}(\oplus_{i=0}^n\pi\widetilde{\mathcal{U}}(-i),\pi M)\\
&=& \bigoplus_{0= i}^n\mathbf{R}Hom_{qgr(\widetilde{\mathcal{U}})}(\pi\widetilde{\mathcal{U}}(-i),\pi M)\\
&=& \bigoplus_{0= i}^n \mathbf{R}Hom_{gr(\widetilde{\mathcal{U}})}(\widetilde{\mathcal{U}},\mathbf{R}\omega\pi M (i))\\
&=& \bigoplus_{0= i}^n
[\mathbf{R}\omega\pi(M)]_{-i}
\end{eqnarray*}
The extra structure needed to make $\bigoplus_{i=0}^n[\mathbf{R}\omega\pi(M)]_{-i}$ into a derived left $E$-module is the collection of action maps
\[\widetilde{\mathcal{U}}^{j-i}\otimes_X [\mathbf{R}\omega\pi(M)]_{-j}\rightarrow [\mathbf{R}\omega\pi(M)]_{-i}\]
which come from $\mathbf{R}\omega\pi(M)$'s left $\widetilde{\mathcal{U}}$-module structure.

Either way one writes it, it defines an equivalence of derived categories.

\begin{mthm}
The functor $\mathbf{R}Hom_{qgr}(T,-)=\bigoplus_{0= i}^n
[\mathbf{R}\omega\pi(-)]_{-i}$ defines an equivalence of triangulated categories (in fact, of dg categories)
\[D^b(qgr(\widetilde{\mathcal{U}}))\simeq D^b(mod(E))\]
with inverse given by $T\otimes_E^\mathbf{L} -$.
\end{mthm}
\begin{proof}
The theorem will follow from the following lemma.
\begin{lemma}
Let $\mathcal{A}$ be an abelian category, and let $T$ be an object in $\mathcal{A}$ which is:
\begin{itemize}
\item \textbf{Compact}: The functor $Hom_{\mathcal{A}}(T,-)$ commutes with direct sums.
\item \textbf{Generator}: For any object $M\in\mathcal{A}$, there is a surjection $T^{\oplus I}\rightarrow \mathcal{A}$ for some index set $I$.
\item \textbf{Finite Dimension}: There is some $i$ such that $Ext^j_{\mathcal{A}}(T,M)=0$ for all $j>i$ and $M\in \mathcal{A}$.
\item $Ext^i_{\mathcal{A}}(T,T)=0$ for $i>0$.
\end{itemize}
Then $\mathbf{R}Hom_{\mathcal{A}}(T,-)$ defines a quasi-equivalence of dg categories (and hence an equivalence of triangulated categories)
\[D^b(\mathcal{A})\simeq D^b(mod(End(T)^{op}))\]
with inverse $T\otimes_{End(T)^{op}}^{\mathbf{L}} -$.
\end{lemma}
\begin{proof}
Theorem 4.3  in \cite{KellerDGcategories} (see also Theorem 8.7 in \cite{KellerDGTilting}) provides a a quasi-equivalence of dg categories $D^b(\mathcal{A})\simeq Perf(Mod(End(T)^{op}))$, where $Perf(Mod(E))$ is the category of perfect complexes.  However, by the finite dimensionality, the image of the functor takes bounded complexes to bounded complexes.  Therefore, $Perf(Mod(End(T)^{op}))\simeq D^b(mod(End(T)^{op}))$.
\end{proof}

The compactness of $T$ is immediate, because $\pi$ is a compact functor and $T$ is $\pi$ of a f.g. object.  The fact that $T$ is a generator was Corollary \ref{coro:gen T}.  The Beilinson resolution proves that $\mathbf{R}\omega\pi$ has finite homological dimension (though it does not give a sharp bound), and so then $\mathbf{R}Hom_{qgr}(T,-)$ does as well.  Finally, the vanishing of higher $Ext$s was Theorem \ref{thm:end T}.
\end{proof}

One interpretation of this theorem is that an object $\pi M\in qgr(\widetilde{\mathcal{U}})$ can be completely determined by knowing $\mathbf{R}\omega\pi(M)$ in degrees $-n$ to $0$, together with knowing the action maps
\[\widetilde{\mathcal{U}}^{j-i}\otimes_X [\mathbf{R}\omega\pi(M)]_{-j}\rightarrow [\mathbf{R}\omega\pi(M)]_{-i}\]
In fact, any object in $D^b(qgr(\widetilde{\mathcal{U}}))$ can be constructed by giving $n+1$ objects $N_{-i}\in D^b(\mathcal{O}_X)$, together with action maps $\widetilde{\mathcal{U}}^{j-i}\otimes_X N_{-j}\rightarrow N_{-i}$ which are required to be associative in the natural way.

This can even be an effective method for constructing objects in $qgr(\widetilde{\mathcal{U}})$, provided one has some method of ensuring that the higher cohomologies vanish.  This was the method used in \cite{BZN} to construct the moduli space of left ideals in the ring of differential operators on a curve.

\section{Examples and Applications.}\label{section:examples}

The generality of Lie algebroids means that this theorem encompasses a wide array of different examples.  We review some of these examples now.

\subsection{Example: Polynomial Algebra.}

This is the case $X=Spec(\mathbb{C})$, and $L$ abelian.  $L$ is then a vector space with trivial Lie bracket.  If $\{x_1,x_2,...x_n\}$ is a basis for $L$, $\mathcal{U}$ is $\mathbb{C}[x_1,x_2,...,x_n]$ and $\widetilde{\mathcal{U}}=\mathbb{C}[t,x_1,x_2,...x_n]$.  Therefore, $qgr(\widetilde{\mathcal{U}})=mod(\mathbb{P}^n)$, by the projective Serre equivalence.  Then the main theorem becomes the derived equivalence of $\mathbb{P}^n$ and the algebra
\[\left(\begin{array}{ccccc}
\mathbb{C} & \mathbb{C}\oplus L & \mathbb{C}\oplus L\oplus Sym^2L & \cdots & \mathbb{C}\oplus L\oplus ...\oplus Sym^nL \\
0 & \mathbb{C} & \mathbb{C}\oplus L & \cdots & \mathbb{C} \oplus L\oplus ... \oplus Sym^{n-1}L \\
0 & 0 & \mathbb{C} & \cdots & \mathbb{C} \oplus L \oplus ... \oplus Sym^{n-2}L \\
\vdots & \vdots & \vdots & \ddots &\vdots \\
0& 0& 0& \cdots & \mathbb{C} \\
\end{array}\right)\]
This algebra is usually written as the path algebra of a quiver $Q_n$, called the \textbf{$n$th Beilinson quiver}.  The equivalence $D^b(mod(\mathbb{P}^n))\simeq D^b(mod(Q_n))$ is the original Beilinson equivalence, and was proven in the seminal paper \cite{BeilinsonEquivalence}.

\subsection{Example: Enveloping Algebra of a Lie Algebra.}

This is the case $X=Spec(\mathbb{C})$, and $L=\mathfrak{g}$, some Lie algebra.  The enveloping algebra is then the usual enveloping algebra $\mathcal{U}\mathfrak{g}$ of the Lie algebra, and $\widetilde{\mathcal{U}\mathfrak{g}}$ is the homogenization.  The categories $qgr(\widetilde{\mathcal{U}\mathfrak{g}})$ were first introduced by \cite{LeBruynvdB} under the name \textbf{quantum space of a Lie algebra}.  The main theorem becomes the derived equivalence of this category and the algebra
\[\left(\begin{array}{ccccc}
\mathbb{C} & (\mathcal{U}\mathfrak{g})^1 & (\mathcal{U}\mathfrak{g})^2 & \cdots & (\mathcal{U}\mathfrak{g})^n \\
0 & \mathbb{C} & (\mathcal{U}\mathfrak{g})^1 & \cdots & (\mathcal{U}\mathfrak{g})^{n-1} \\
0 & 0 & \mathbb{C} & \cdots & (\mathcal{U}\mathfrak{g})^{n-2} \\
\vdots & \vdots & \vdots & \ddots &\vdots \\
0& 0& 0& \cdots & \mathbb{C} \\
\end{array}\right)\]
This algebra again can be written as the path algebra of a quiver, which will look like the $n$th Beilinson quiver with its relations deformed by the Lie bracket.

\subsection{Example: Differential Operators.}

In this case, $X$ is any irreducible smooth affine variety, and $L$ is the tangent bundle $\mathcal{T}$.  Then, $\mathcal{U}$ is $\mathcal{D}$, the ring of differential operators, and $\widetilde{\mathcal{U}}$ is $\widetilde{\mathcal{D}}$, the Rees algebra of the differential operators.  The category $qgr(\widetilde{\mathcal{D}})$ is then derived equivalent to the algebra
\[\left(\begin{array}{ccccc}
\mathcal{O}_X & \mathcal{D}^1 & \mathcal{D}^2 & \cdots & \mathcal{D}^d \\
0 & \mathcal{O}_X & \mathcal{D}^1 & \cdots & \mathcal{D}^{d-1} \\
0 & 0 & \mathcal{O}_X & \cdots & \mathcal{D}^{d-2} \\
\vdots & \vdots & \vdots & \ddots & \vdots \\
0 & 0 & 0 & \cdots & \mathcal{O}_X \\
\end{array}\right)\]
Not much else can be said in this level of generality.  However, for a powerful application of this in the form of curves, see Subsection \ref{subsection:ideals}.

\subsection{Non-Examples.}\label{subsection:nonexamples}

It is worth noting that $\widetilde{\mathcal{U}}$ is not the most general class of graded algebra for which the techniques here work, and for which a similar version of the main theorem applies.  For example, let $PP_\hbar$ denote the algebra over $\mathbb{C}$ generated by $w_1$, $w_2$, and $w_3$, subject to the relations
\[ [w_1,w_3]=[w_2,w_3]=0,\;\;\; [w_1,w_2]=2\hbar w_3^2\]
One can check that this is not the homogenization of any universal enveloping algebra of a Lie algebra.

However, in \cite{KKO}, a similar Koszul theory is developed, as well as a similar Beilinson equivalence, which is then used for a monad-theoretic construction of the moduli space of certain kinds of modules.

Another non-example of a relatively quadratic algebra which has an identical Koszul theory and Beilinson transform is the $\widetilde{\mathcal{U}}^{op}$, the opposite algebra of the enveloping algebra of a Lie algebroid.  This is equivalent to showing that the category of graded {right} $\widetilde{\mathcal{U}}$-module has a quotient $qgr(\widetilde{\mathcal{U}})$ which satisfies all the theorems of this paper.  Every proof in this paper works in this case, occasionally with slight modification (actually, the proof of the relative Gorenstein property is a little bit shorter).

So then, \textbf{what is the most general setting where the proofs in this paper work?}  The answer is that the proofs in this paper will work for any relatively quadratic algebra $A$, such that
\begin{itemize}
\item $A$ is \textbf{Koszul}, in that the left and right Koszul complexes are resolutions of $\mathcal{O}_X$.
\item $A^\perp$ is a finitely generated projective left and right $\mathcal{O}_X$-module and relatively Frobenius over $\mathcal{O}_X$.  That is, Corollaries \ref{coro:dual fg proj}, \ref{coro:dual finite} and \ref{coro:dual duality} hold.
\end{itemize}
The results can be generalized in a different direction, by extending algebras to sheaves of algebras on a non-affine $X$.  See Appendix \ref{section:local} for details.

\subsection{Application: Grothendieck Group.}

An immediate application of the derived equivalence is computing the Grothendieck group $K(qgr(\widetilde{\mathcal{U}}))$ of the category $qgr(\widetilde{\mathcal{U}})$, because the Grothendieck group depends only on the bounded derived category.  Furthermore, $K(mod(E))$ is easy to compute because, like a quiver, it can be shown that the Grothendieck group depends only on the diagonal part of $E$ (the vertices) and not on the above diagonal part (the arrows).

\begin{lemma}
$K(mod(E))=K(coh(X))^{\oplus (n+1)}$.
\end{lemma}
\begin{proof}
Let $M\in mod(E)$, and let $e_{-i}$ denote the idempotent in $E$ which is $1\in\mathcal{O}_X$ in the $(n+1-i,n+1-i)$ entry in the matrix.  Recall that $M$ can be described by the $\mathcal{O}_X$-modules $M_{-i}:=e_iM\in coh(X)$, together with a collection of action maps $\mathcal{U}^1\otimes_X M_{-i}\rightarrow M_{-i+1}$.  Note that $M$ has a filtration by submodules $M_{\geq -i}:= (\sum_{j=0}^ie_{-i})M$, with the action maps the same as $M$ where they aren't necessarily zero.  The successive quotients $M_{\geq -i}/M_{\geq -i+1}= M_{-i}$, and so $[M]=\sum_{i=0}^n[M_{-i}]$.  Therefore, $K(mod(E))$ is generated by the class of modules of the form $M_{-i}$ for some $M$.

Let $N$ and $N'$ be two $\mathcal{O}_X$ modules, and let $e_{-i}N$ and $e_{-i}N'$ be the corresponding $E$-modules.  Then $[e_{-i}N]=[e_{-i}N']$ only if $[N]=[N']$ in $K(coh(X))$.  Furthermore, $[e_{-i}N]=[e_{-j}N']$ for $i\neq j$ only if both are the zero class.  Therefore, the group $K(mod(E))$ decomposes into $K(coh(X))^{\oplus (n+1)}$, where $[M]$ goes to $([M_{0}],[M_{-1}],...[M_{-n}])$.
\end{proof}

\begin{coro}
$K(qgr(\widetilde{\mathcal{U}}))\simeq K(coh(X))^{\oplus(n+1)}$.
\end{coro}

Explicitly, under this isomorphism, $[\pi M]$ goes to 
\[([\mathbf{R}\omega\pi(M)_0],[\mathbf{R}\omega\pi(M)_{-1}],... [\mathbf{R}\omega\pi(M)_{-n}])\]
This decomposition can be used to define the notion of a $K(coh(X))$-valued $i$th Chern class for an object in $qgr(\widetilde{\mathcal{U}})$.  Let the $i$-th Chern class of $\pi M$ be defined as
\[ \sum_{j=0}^n { i \choose j} [\mathbf{R}\omega\pi(M)_{-j}] \in K(coh(X))\]
where ${i\choose j}=0$ if $j>i$.  In the case of $\mathbb{P}^n$, this will coincide with the usual Chern class of a module.  

\subsection{Application: Ideals in $\mathcal{U}$.}\label{subsection:ideals}

So far, the study of the category $qgr(\widetilde{\mathcal{U}})$ has been motivated by its appealing properties (the main theorem, for instance), and by its close but nebulous relation with the study of $mod(\mathcal{U})$.  We now briefly illustrate an example of the latter, by showing how left ideals in $\mathcal{U}$ can be studied via this method.

A left ideal $I\subset \mathcal{U}$ comes naturally equipped with a filtration, as a restriction of the filtration on $\mathcal{U}$.  This translates into a short exact sequence in $gr(\widetilde{\mathcal{U}})$:
\[0\rightarrow \widetilde{\mathcal{I}}\rightarrow \widetilde{\mathcal{U}}\rightarrow \widetilde{\mathcal{U}/I}\rightarrow 0\]
We then wish to study $I$ by studying $\pi\widetilde{I}$, but we need to be able to recover $I$ from $\pi\widetilde{I}$.
\begin{lemma}
$\omega\pi \widetilde{I}\simeq \widetilde{I}$.
\end{lemma}
\begin{proof}
By the exact sequence
\[0\rightarrow \tau(\widetilde{I})\rightarrow \widetilde{I}\rightarrow \omega\pi \widetilde{I} \rightarrow \mathbf{R}^1\tau \widetilde{I}\rightarrow 0,\]
it suffices to show that $\tau \widetilde{I}=\mathbf{R}^1\tau(\widetilde{I})=0$.  Note that for any filtered $\mathcal{U}$-module $M$, $t$ acts as an inclusion on $\widetilde{M}$.  Therefore, $\underline{Hom}_{gr(\widetilde{\mathcal{U}})}(\widetilde{\mathcal{U}}/\widetilde{\mathcal{U}}^{\geq n},\widetilde{M})=0$, and so $\tau \widetilde{M}=0$.  This means that $\tau\widetilde{I}=0$.

Now, apply $\mathbf{R}\tau$ to the sequence
\[0\rightarrow \widetilde{I}\rightarrow \widetilde{\mathcal{U}}\rightarrow \widetilde{\mathcal{U}/I}\rightarrow 0\]
Since $\tau\widetilde{\mathcal{U}}=\mathbf{R}^1\tau(\widetilde{\mathcal{U}})=0$ by Lemma \ref{lemma: tau vanishing}, we know that $\mathbf{R}^1\tau(\widetilde{I})=\tau(\widetilde{\mathcal{U}/I})$, but $\tau$ of anything which is the homogenization of a filtered module is zero.  Therefore, $\mathbf{R}^1\tau (\widetilde{I})=0$.
\end{proof}

This was the approach used by \cite{BZN} to characterize ideal classes in the ring of differential operators on a curve $X$.  The general idea is to characterize which derived $E$-modules came from ideal classes, and show that every such derived $E$ module came from an ideal class.
\begin{thm}(\cite{BZN}, Theorem 4.3)
Let $I$ be an ideal in $\mathcal{D}$ for $X$ a smooth curve.  Then
\begin{enumerate}
\item $(\mathbf{R}\omega\pi(\widetilde{I}))_{-1}=V [-1]$, where $V$ is a finite-length sheaf on $X$.
\item $(\mathbf{R}\omega\pi(\widetilde{I}))_0=Cone(i:J\rightarrow V)$, where $J$ is some ideal on $X$ and $i$ is some $\mathcal{O}_X$-module map.
\item The action map $a:\mathcal{D}^1\otimes_X (\mathbf{R}\omega\pi(\widetilde{I}))_{-1}\rightarrow (\mathbf{R}\omega\pi(\widetilde{I}))_0$ restricts on $\mathcal{O}_X$ to the natural map
    \[\left\{\begin{array}{c} V \\ \uparrow \\ 0 \end{array}\right\}\begin{array}{c} \rightarrow_{Id_V} \\ \\ \rightarrow_0\end{array} \left\{\begin{array}{c} V \\ \uparrow \\ J \end{array}\right\}\]
\end{enumerate}
Furthermore, any choice of such $V$, $J$, $i$ and $a$ will determine a derived $E$-module which corresponds to an ideal under the inverse Beilinson equivalence.
\end{thm}

\appendix

\section{The Quadratic Dual.}
This appendix collects and proves the important facts about the quadratic dual algebra, $\widetilde{\mathcal{U}}^\perp$.
\subsection{The Structure of the Quadratic Dual.}
This section explores the structure of $\widetilde{\mathcal{U}}^{\perp}$ as an algebra.  First, note that $J^1$ fits into a short exact sequence of $\mathcal{O}_X$-bimodules,
\[0\rightarrow L^*\rightarrow J^1\rightarrow \mathcal{O}_X\rightarrow 0\]
The 'action on 1' map $\mathcal{U}\rightarrow \mathcal{O}_X$ is a map of left $\mathcal{U}$-modules.  It restricts to a map of left $\mathcal{O}_X$-modules $e:\mathcal{U}^1\rightarrow \mathcal{O}_X$, and so it determines an element $e\in J^1$ and its image in $\widetilde{\mathcal{U}}^{\perp}$.  Since $e$ acts as the identity on $\mathcal{O}_X\subset \mathcal{U}^1$, its image under the map $J^1\rightarrow \mathcal{O}_X$ is the identity in $\mathcal{O}_X$.

Next, define the \textbf{$L$-exterior derivative} $\mu:L^*\rightarrow  L^*\otimes_X L^*=(L\otimes_XL)^*$ by 
\[\mu(\sigma)(l\otimes l'):= \frac{1}{2}\left[d_{\tau(l)}(\sigma(l')) -d_{\tau(l')}(\sigma(l))-\sigma([l,l'])\right]\]
The name comes from the case when $L=\mathcal{T}$, where $\mu:\mathcal{T}^*\rightarrow\mathcal{T}^*\otimes_X \mathcal{T}^*$ is the usual exterior derivative.

Finally, a quick rundown on the explicit form of some definitions for those who haven't had the luxury to work out examples.
\begin{itemize}
\item The way the $\mathcal{O}_X$-bimodule structure on $J^1=\,^*(\mathcal{U}^1)$ was defined, $(ae)(\partial)=e(\partial a)$.  
\item From the isomorphism $\,^*(\mathcal{U}^1)\otimes_X \,^*(\mathcal{U}^1)=\,^*(\mathcal{U}^1\otimes_X \mathcal{U}^1)$, for $\sigma,\sigma'\in \,^*(\mathcal{U}^1)$, $(\sigma\otimes \sigma')(\partial\otimes \partial')=\sigma'(\partial\cdot\sigma(\partial'))$.
\item From the definition of $\mathcal{U}$, we see that for $\partial\in ker(e)$ and $a\in \mathcal{O}_X$, then $[\partial,a]=d_{\tau(\partial)}(a)$.
\end{itemize}

\begin{lemma}The element $e\in \widetilde{\mathcal{U}}^{\perp}$ satisfies
\begin{enumerate}
\item $e^2=0$.
\item $ae-ea=\tau^\vee(da)$, for $a\in\mathcal{O}_X$, and where $\tau^\vee:\mathcal{T}^*\rightarrow L^*$ is dual to the anchor map $L\rightarrow \mathcal{T}$.
\item $\sigma e+e\sigma=\mu({\sigma})$, for $\sigma\in L^*\in J^1$.
\end{enumerate}
\end{lemma}
\begin{proof}
The easy relation to show is $(2)$, because it is a degree 1 relation.  Consider the element $ae -ea\in J^1$, and apply it to any $\partial\in \mathcal{U}^1$.
\[(ae-ea)\partial=e(\partial a)-e(a\partial )=e([\partial,a])=d_{\tau(\partial)}(a)=\iota_{da}(\tau(\partial))=\tau^\vee(da)\partial\]
and so $(ae-ea)=\tau^\vee(da)$.

The other two relations are degree 2, so they are true if and only if they are in $J^2$; that is, if they kill $R\in \mathcal{U}^1\otimes_X \mathcal{U}^1$.  Remember that $R$ is spanned by elements of the form $\partial\otimes \partial' -\partial'\otimes \partial -[\partial,\partial']\otimes 1$.

\underline{(1) $e\otimes e$}.
\[(e\otimes e)(\partial\otimes \partial' -\partial'\otimes \partial-[\partial,\partial']\otimes 1)\]
\begin{eqnarray*}&=&e(\partial e(\partial')) - e(\partial' e(\partial)) - e([\partial,\partial'] e(1))\\
&=&e(\partial')e(\partial) + e([\partial,e(\partial')])-e(\partial)e(\partial') - e([\partial',e(\partial)]) -e([\partial,\partial'])\\
&=&[\partial,e(\partial')]-[\partial',e(\partial)] -e([\partial,\partial'])
\end{eqnarray*}
It suffices to check that this final expression vanishes in several cases.  
\begin{itemize}
\item If both $\partial$ and $\partial'$ are in $\mathcal{O}_X$, then all the commutators vanish.  
\item If one of $\partial$ and $\partial'$ is in $\mathcal{O}_X$ and the other is in the kernel of $e$, then one of the terms vanish and the other two terms are identical.  
\item If both $\partial$ and $\partial'$ are in the kernel of $e$, then this is also true of their commutator, and so all three terms vanish.
\end{itemize}

\underline{(3) ${\sigma}\otimes e +e\otimes {\sigma} - \mu({\sigma})$}.
\[({\sigma}\otimes e +e\otimes \sigma)(\partial\otimes \partial' -\partial'\otimes \partial-[\partial,\partial']\otimes 1)\]
\begin{eqnarray*}
&=&\left[e(\partial {\sigma}(\partial'))- e(\partial' {\sigma}(\partial)) \right]+\left[{\sigma}(\partial e(\partial'))-{\sigma}(\partial' e(\partial)) -{\sigma}([\partial,\partial'] e(1))\right] \\
&=&e(\partial{\sigma}(\partial'))-e(\partial'{\sigma}(\partial))+e(\partial'){\sigma}(\partial)-e(\partial){\sigma}(\partial')-{\sigma}([\partial,\partial'])\\
&=&e([\partial',{\sigma}(\partial)])-e([\partial,{\sigma}(\partial')])-{\sigma}([\partial,\partial'])\\
&=&[\partial',{\sigma}(\partial)]-[\partial,{\sigma}(\partial')]-{\sigma}([\partial,\partial'])\\
&=&d_{\tau(\partial)}(\sigma(\partial')) - d_{\tau(\partial')}(\sigma(\partial)) - \sigma([\partial,\partial'])
\end{eqnarray*}
Compare to 
\[\mu(\sigma)(\partial\otimes \partial' -\partial'\otimes \partial-[\partial,\partial']\otimes 1)\]
\begin{eqnarray*}
&=& \frac{1}{2}\left[d_{\tau(\partial)}(\sigma(\partial'))-d_{\tau(\partial')}(\sigma(\partial)) - \sigma([\partial,\partial'])\right] \\
&-& \frac{1}{2}\left[d_{\tau(\partial')}(\sigma(\partial)) +d_{\tau(\partial)}(\sigma(\partial')) +\sigma([\partial',\partial]) \right] \\
&-& \frac{1}{2}\left[d_{\tau([\partial,\partial'])}(\sigma(1)) - d_{\tau(1)}([\partial,\partial']) - \sigma( [[\partial,\partial'],1]) \right]\\
&=&d_{\tau(\partial)}(\sigma(\partial')) - d_{\tau(\partial')}(\sigma(\partial)) - \sigma([\partial,\partial']) 
\end{eqnarray*}
Therefore, $\sigma\otimes e+e\otimes \sigma-\mu(\sigma)$ kills $R\in \mathcal{U}^1\otimes_X\mathcal{U}^1$, and so it is a relation in $\widetilde{\mathcal{U}}^{\perp}$.
\end{proof}

For any element $\widetilde{\mathcal{U}}^{\perp}$, the above (graded) commutators allow $e$ to collected on one side (for instance, to the right).  Since $e^2=0$, an element in $\widetilde{\mathcal{U}}^{\perp}$ can have at most one $e$ in it.  The following theorem then establishes that $\widetilde{\mathcal{U}}^{\perp}$ is a rank 2 module over the subalgebra of elements without an $e$.

\begin{thm}
The map $L^* \rightarrow J^1$ extends to an inclusion $\Lambda_X^\bullet L^*\rightarrow \widetilde{\mathcal{U}}^{\perp}$.  This map fits into a short exact sequence of graded $\Lambda_X^\bullet L^*$-bimodules
\[0\rightarrow \Lambda_X^\bullet L^* \rightarrow \widetilde{\mathcal{U}}^{\perp} \rightarrow \Lambda_X^\bullet L^*(-1)\rightarrow 0\]
where $e\in \widetilde{\mathcal{U}}^{\perp}$ goes to $1\in \Lambda_X^\bullet L^* (-1)$.
\end{thm}
\begin{proof}
First, it is easy to see that, for $\sigma,\sigma'\in L^*$, $\sigma\otimes \sigma'+\sigma'\otimes\sigma$ is a relation in $\widetilde{\mathcal{U}}^{\perp}$.
\[(\sigma\otimes \sigma'+\sigma'\otimes\sigma)(\partial\otimes\partial'-\partial'\otimes \partial -[\partial,\partial']\otimes1)\]
\begin{eqnarray*}
&=& \sigma'(\partial \sigma(\partial')) +\sigma(\partial\sigma'(\partial')) -\sigma'(\partial'\sigma(\partial)) -\sigma(\partial\sigma'(\partial'))\\
&=& \sigma(\partial')\sigma'(\partial) + \sigma'(\partial')\sigma(\partial) -\sigma(\partial)\sigma'(\partial') -\sigma'(\partial')\sigma(\partial) =0
\end{eqnarray*}
It is not much harder to see that any relation in $L^*\otimes_X L^*$ is fixed by the map which sends $\sigma\otimes \sigma'$ to $\sigma'\otimes \sigma$.  Therefore, elements of the form $\sigma\otimes\sigma'+\sigma'\otimes \sigma$ generate the relations in $L^*\otimes_X L^*$.  It follows that the submodule $L^*\subset J^1\subset \widetilde{\mathcal{U}}^{\perp}$ generates a copy of the algebra $\Lambda_X^\bullet L^*$.

Now, let $C$ denote the cokernel of $\Lambda_X^\bullet L^*\rightarrow \widetilde{\mathcal{U}}^{\perp}$, as a $\Lambda^\bullet_X L^*$-bimodule.  Note that the previous lemma showed that the (graded) commutator of $e$ with any element of $J^1$ lies in $L^*\subset J^1$.  Therefore, the image of $e$ in $\widetilde{\mathcal{U}}^{\perp}\rightarrow C$ is (graded) central.  Furthermore, since $e^2=0$, $e$ generates $C$, and so there is a surjective map $\Lambda_X^\bullet L^*(-1)\rightarrow C$ which sends $1$ to $e$.

For this not to be an isomorphism, there would have to be a relation of the form $\sigma e  - \Upsilon$, for $\sigma\in L^*$ and $\Upsilon\in L^*\otimes_X L^*$.  Let $\partial$ be an element in $\mathcal{U}^1$ which is not killed by $\sigma$.  Then
\[ (\sigma\otimes e)(\partial\otimes 1 - 1\otimes \partial)=e(\partial \sigma(1))-e(\sigma(\partial))=\sigma(\partial) \]
By construction, this is not zero.  However, $\Upsilon$ must kill $\partial\otimes 1-1\otimes \partial$ since $L^*$ kills $1\in \mathcal{U}^1$.  Therefore, there cannot be such a relation, and the map $\Lambda^\bullet_X L^*(-1)\rightarrow C$ is an isomorphism.
\end{proof}

Since $\Lambda_X^\bullet L^*$ is an algebra which is finitely generated projective as a $\mathcal{O}_X$-module on either side and zero in large enough degree, we can deduce identical facts about $\widetilde{\mathcal{U}}^{\perp}$.

\begin{coro}\label{coro:dual fg proj}
For all $i$, $\widetilde{\mathcal{U}}^{\perp i}$ is a finitely generated, projective $\mathcal{O}_X$-module on the left and right.
\end{coro}

\begin{coro}\label{coro:dual finite}
If $i>n+1$, then $\widetilde{\mathcal{U}}^{\perp i}=0$.
\end{coro}

Let $\omega_L$ denote $\Lambda^n_XL^*$, the top exterior power of the dual to $L$.  From the Lemma, it is clear that $\widetilde{\mathcal{U}}^{\perp n+1}=\omega_L$.  This now gives a pairing between elements of $\widetilde{\mathcal{U}}^{\perp}$ whose degree adds to $n+1$.  We then have
\begin{lemma} (The relative Frobenius property)
For any $i$, the multiplication map
\[\widetilde{\mathcal{U}}^{\perp i}\otimes_X \widetilde{\mathcal{U}}^{\perp (n+1-i)}\rightarrow \omega_L\]
is a 'perfect pairing'.  That is, the adjoint maps
\[\widetilde{\mathcal{U}}^{\perp (n+1-i)}\rightarrow Hom_{X-}(\widetilde{\mathcal{U}}^{\perp i},\omega_L),\;\;\; and \;\;\;\widetilde{\mathcal{U}}^{\perp i}\rightarrow Hom_{-X}(\widetilde{\mathcal{U}}^{\perp (n+1-i)},\omega_L)\] are isomorphisms of $\mathcal{O}_X$-bimodules.
\end{lemma}
\begin{proof}
Explicitly, the adjoint map $\widetilde{\mathcal{U}}^{\perp (n+1-i)}\rightarrow Hom_{X-}(\widetilde{\mathcal{U}}^{\perp i},\omega_L)$ takes an element $\mu\in \widetilde{\mathcal{U}}^{\perp (n+1-i)}$ and sends it to the map $\gamma\in\widetilde{\mathcal{U}}^{\perp i}\rightarrow \mu\cdot\gamma\in \omega_L$. Consider the short exact sequence of $\mathcal{O}_X$-bimodules
\[0\rightarrow \Lambda^{n+1-i}_XL^*\rightarrow \widetilde{\mathcal{U}}^{\perp n+1-i}\rightarrow \Lambda^{n+i}_XL^*\rightarrow 0\]
If $\mu\in \Lambda^{(n+1-i)}_XL^*\in\widetilde{\mathcal{U}}^{\perp (n+1-i)}$, then $\mu\cdot \gamma$ only depends on the image of $\gamma$ under the map $\widetilde{\mathcal{U}}^{\perp i}\rightarrow \Lambda^{i-1}_XL^*$.  Similarly, if we know that $\gamma \in \Lambda_X^{i}L^*\subset\widetilde{\mathcal{U}}^{\perp i}$, then $\mu\cdot\gamma$ only depends on the image of $\mu$ under the map  $\widetilde{\mathcal{U}}^{\perp (n+1-i)}\rightarrow \Lambda^{n+i}_XL^*$.  This means that the adjoint map above splits into a map of short exact sequences
\[\begin{array}{ccccc}
\Lambda^{n+1-i}_XL^* & \rightarrow & \widetilde{\mathcal{U}}^{\perp (n+1-i)} & \rightarrow & \Lambda^{n-i}_XL^*  \\
\downarrow & & \downarrow & & \downarrow \\
Hom_{X-}(\Lambda^{i-1}_XL^*,\omega_L) & \rightarrow & Hom_{X-}(\widetilde{\mathcal{U}}^{\perp i},\omega_L) & \rightarrow & Hom_{X-}(\Lambda^{i}_XL^*,\omega_L)
\end{array}\]
The left and right maps are isomorphisms, because they are both adjoint to multiplication maps of the form $\Lambda^j_XL^*\otimes_X \Lambda^{n-j}_XL^*\rightarrow \omega_L$.  Therefore, the middle map is an isomorphism.  The proof for the other adjoint map is identical.
\end{proof}

This can be restated in a more compact form.
\begin{coro}\label{coro:dual duality}
There are isomorphisms of $\mathcal{O}_X$-bimodules
\[(\widetilde{\mathcal{U}}^{\perp i})^*=\omega_L^*\otimes_X  \widetilde{\mathcal{U}}^{\perp (n+1-i)},\;\;\;
\,^*(\widetilde{\mathcal{U}}^{\perp i})=  \widetilde{\mathcal{U}}^{\perp (n+1-i)}\otimes_X\omega_L^*\]
\end{coro}
\begin{proof}
\[\widetilde{\mathcal{U}}^{\perp (n+1-i)}\simeq Hom_{X-}(\widetilde{\mathcal{U}}^{\perp i},\omega_L)=Hom_{X-}(\widetilde{\mathcal{U}}^{\perp i},\mathcal{O}_X)\otimes_X \omega_L=\,^*(\widetilde{\mathcal{U}}^{\perp i})\otimes_X \omega_L\]
Similarly, $\widetilde{\mathcal{U}}^{\perp (n+1-i)}=\omega_L^*\otimes_X (\widetilde{\mathcal{U}}^{\perp})^*$.
Since $\omega_L$ is a line bundle, tensoring these with $\omega_L^*$ on the left or right gives the theorem.
\end{proof}

\subsection{The Quadratic Dual as an $Ext$ Algebra.}

In this section, we prove the following theorem about $\widetilde{\mathcal{U}}^\perp$.  
\begin{thm}\label{Thm:Ext Algebra}
$\widetilde{\mathcal{U}}^\perp$ is isomorphic to $\underline{Ext}^\bullet_{\widetilde{\mathcal{U}}-}(\mathcal{O}_X,\mathcal{O}_X)$ as a graded algebra, where $J^1=\,^*(\mathcal{U}^1)\subset \widetilde{\mathcal{U}}^\perp$ has graded degree -1.
\end{thm}
\begin{proof}
It is easy to see this isomorphism, on the level of graded $\mathcal{O}_X$-modules.
\begin{lemma}\label{lemma:Mod Ext Algebra}
$\widetilde{\mathcal{U}}^\perp$ is isomorphic to $\underline{Ext}^\bullet_{\widetilde{\mathcal{U}}-}(\mathcal{O}_X,\mathcal{O}_X)$ as a graded $\mathcal{O}_X$-module.
\end{lemma}
\begin{proof}
The left Koszul resolution $K^\bullet_{right}$ is a left projective resolution of $\mathcal{O}_X$.  Therefore,
\begin{eqnarray*}
\mathbf{R}\underline{Hom}^\bullet_{\widetilde{\mathcal{U}}-}(\mathcal{O}_X,\mathcal{O}_X)
&=&\underline{Hom}_{\widetilde{\mathcal{U}}-}(K^\bullet,\mathcal{O}_X)\\
&=&\bigoplus_{i=0}^n\underline{Hom}_{\widetilde{\mathcal{U}}-}(\widetilde{\mathcal{U}}(-i) \otimes_X (\widetilde{\mathcal{U}}^{\perp i})^*,\mathcal{O}_X)\\
&=&\bigoplus_{i=0}^nHom_{X-}((\widetilde{\mathcal{U}}^{\perp i})^*,\mathcal{O}_X)(i)\\
&=&\bigoplus_{i=0}^n\widetilde{\mathcal{U}}^{\perp i}(i)
\end{eqnarray*}
Since each term in the complex is concentrated in a different graded degree, the boundary vanishes, and so the cohomology is isomorphic to $\widetilde{\mathcal{U}}^{\perp}$.
\end{proof}

Showing that this is an isomorphism of algebras will require more work.  Let $\mathbf{B}^\bullet$ denote the \textbf{normalized left bar resolution} of $\mathcal{O}_X$ (see \cite{Weibel}, page 284 for details).  This is the complex of graded left $\widetilde{\mathcal{U}}$-modules with $\mathbf{B}^{-i}= \widetilde{\mathcal{U}}\otimes_X(\widetilde{\mathcal{U}}^{\geq 1})^{\otimes_X (i-1)}$
where the boundary sends $a_1\otimes_X a_2 \otimes_X ...\otimes_X a_n$ to
\[ \sum_{i=1}^{n-1}(-1)^ia_1\otimes_Xa_2 \otimes_X ...\otimes_X a_ia_{i+1}\otimes_X ...\otimes_X a_n\]
The complex $\mathbf{B}^\bullet$ is a left projective resolution of $\mathcal{O}_X$, with the augmentation map $\mathbf{B}^0=\widetilde{\mathcal{U}}\rightarrow
\mathcal{O}_X$ the natural projection onto graded degree zero.

Therefore, $\underline{Ext}^\bullet_{\widetilde{\mathcal{U}}-}(\mathcal{O}_X,\mathcal{O}_X)$ is the cohomology algebra of the differential graded algebra (dga) $\underline{Hom}_{\widetilde{\mathcal{U}}-}(\mathbf{B}^\bullet,\mathbf{B}^\bullet)$, where the multiplication is the composition of maps.  The augmentation map $\mathbf{B}^\bullet\rightarrow \mathcal{O}_X$ gives a quasi-isomorphism of complexes
$\underline{Hom}_{\widetilde{\mathcal{U}}-}(\mathbf{B}^\bullet,\mathbf{B}^\bullet)\rightarrow \underline{Hom}_{\widetilde{\mathcal{U}}-}(\mathbf{B}^\bullet,\mathcal{O}_X)$. Since
\begin{eqnarray*}
\underline{Hom}_{\widetilde{\mathcal{U}}-}(\mathbf{B}^{-i},\mathcal{O}_X)
&=&\underline{Hom}_{\widetilde{\mathcal{U}}-}( \widetilde{\mathcal{U}}\otimes_X(\widetilde{\mathcal{U}}^{\geq 1})^{\otimes_X (i-1)},\mathcal{O}_X)\\
&=&\underline{Hom}_{X-}( (\widetilde{\mathcal{U}}^{\geq1})^{\otimes_X (i-1)},\mathcal{O}_X)\\
&=&[\,^*(\widetilde{\mathcal{U}}^{\geq 1})]^{\otimes_X (i-1)}
\end{eqnarray*}
Thus, $\underline{Hom}_{\widetilde{\mathcal{U}}-}(\mathbf{B}^\bullet,\mathcal{O}_X)$ is isomorphic to $T_X\,^*(\widetilde{\mathcal{U}}^{\geq 1})$ as a graded $\mathcal{O}_X$-module, and the natural multiplication on the tensor algebra makes it into a dga.

In fact, the quasi-isomorphism \[\underline{Hom}_{\widetilde{\mathcal{U}}-}(\mathbf{B}^\bullet,\mathbf{B}^\bullet)\rightarrow \underline{Hom}_{\widetilde{\mathcal{U}}-}(\mathbf{B}^\bullet,\mathcal{O}_X)=T_X\,^*(\widetilde{\mathcal{U}}^{\geq 1})\] is a map of dgas.  To see this, let us construct a section of this map.  Let $\phi\in[\,^*(\widetilde{\mathcal{U}}^{\geq 1})]^{\otimes_X (i-1)}$, then for any $j>i$, there is a natural map
\[ \widetilde{\mathcal{U}}\otimes_X(\widetilde{\mathcal{U}}^{\geq 1})^{\otimes_X (j-1)}\rightarrow  \widetilde{\mathcal{U}}\otimes_X(\widetilde{\mathcal{U}}^{\geq 1})^{\otimes_X (j-i-1)}\]
given by applying $\phi$ to the first $i$ terms on the left.  It is easy but tedious to verify that this gives a map of dgas $T_X\,^*(\widetilde{\mathcal{U}}^{\geq 1})\rightarrow\underline{Hom}_{\widetilde{\mathcal{U}}-}(\mathbf{B}^\bullet,\mathbf{B}^\bullet)$ which is a section of the above map.  Therefore, $\underline{Ext}^\bullet_{\widetilde{\mathcal{U}}-}(\mathcal{O}_X,\mathcal{O}_X)$ is the cohomology algebra of the dga $T_X\,^*(\widetilde{\mathcal{U}}^{\geq 1})$.

The dga $T_X\,^*(\widetilde{\mathcal{U}}^{\geq 1})$ has both a cohomological degree (coming from the usual grading on a tensor algebra) and a graded degree (coming from the grading on $\widetilde{\mathcal{U}}^{\geq 1}$).  Because $\,^*(\widetilde{\mathcal{U}}^{\geq 1})$ is concentrated in graded degree $\leq -1$, $[\,^*(\widetilde{\mathcal{U}}^{\geq 1})]^{\otimes_X (i-1)}$ is concentrated in graded degree $\leq -i$.  Therefore, if one restricts the complex $T_X\,^*(\widetilde{\mathcal{U}}^{\geq 1})$ to graded degree $-i$, the resulting complex is non-zero in cohomological degrees $j$, $0\leq j\leq i$.

However, we do actually know the cohomology of this complex, due to Lemma \ref{lemma:Mod Ext Algebra}.  Specifically, we know that in graded degree $-i$, the cohomology is concentrated in cohomological degree $i$.  Since the corresponding complex is concentrated in cohomological degrees $\leq i$, the cohomology must be the cokernel of the boundary map.  We therefore have a map of dgas $T_X\,^*(\widetilde{\mathcal{U}}^{\geq1})\rightarrow \widetilde{\mathcal{U}}^{\perp}$, which is a quasi-isomorphism.

Note that, for an element in $T_X\,^*(\widetilde{\mathcal{U}}^{\geq 1})$ to have graded degree $-i$ and cohomological degree $i$, it must be the tensor product of $i$ elements of graded degree $-1$ elements; therefore, $(T_X\,^*(\widetilde{\mathcal{U}}^{\geq1}))^{(-i,i)}=[\,^*(\widetilde{\mathcal{U}}^{1})]^{\otimes_X i}=(J^1)^{\otimes_X i}$.  If we let $T_XJ^1$ be a dga with zero boundary, this extends to a map of dgas $T_XJ^1\rightarrow T_X\,^*(\widetilde{\mathcal{U}}^{\geq1})$, which is the identity in degree $(-i,i)$ and zero elsewhere.

The composition
\[T_XJ^1\rightarrow T_X\,^*(\widetilde{\mathcal{U}}^{\geq1})\rightarrow \widetilde{\mathcal{U}}^\perp\]
is then a surjection of dgas; since their boundaries are zero, we can think of them as algebras again.  Since it is an isomorphism in graded degree $-1$ on the $J^1$'s, its kernel must be exactly generated by $J^2\subset J^1\otimes_X J^1$.  The theorem follows.
\end{proof}

\section{Localization and Sheafification.}\label{section:local}

The results of this paper localize and sheafify correctly, provided one defines the non-affine versions of the constructions correctly.  In this section, we provide the necessary definitions and sketch the necessary proofs.  In this section, $X$ is no longer assumed to be affine, though it is still smooth and irreducible.

\subsection{Lie Algebroids.}

As was mentioned before, Lie algebroids are compatible with localization; that is, the localization of a Lie algebroid naturally has a Lie algebroid structure.  To wit, let $X'$ be an affine open subscheme of affine $X$ defined by a multiplicative subset $S$ of $\mathcal{O}_X$, and let $L':=\mathcal{O}_{X'}\otimes_X L$.
\begin{lemma}
If $(X,L)$ is a Lie algebroid, then $(X',L')$ has a unique Lie algebroid structure which is compatible with the inclusion $L\rightarrow L'$.
\end{lemma}
\begin{proof}
For any $l\in L$ and $s\in S$, the anchor map defines the derivative of $s$ along $l$ to be $d_{\tau(l)}(s)$.  Therefore, there is only one choice for the derivative of $s^{-1}$ along $l$,
\[d_{\tau(l)}(s^{-1}):=-s^{-2}d_{\tau(l)}(s)\]
because $d_{\tau(l)}$ must be a derivation.  In this way, the anchor map $L\rightarrow \mathcal{T}_X$ extends canonically to an anchor map $L\rightarrow \mathcal{T}_{X'}$.  The $\mathcal{O}_{X'}$-module structure on $\mathcal{T}_{X'}$ means that this map extends uniquely to a map $L'\rightarrow \mathcal{T}_{X'}$.

Elements in $L'$ are of the form $s^{-n}\otimes l$, for $s\in S$ and $l\in L$, and so the compatibility of the anchor map with the Lie bracket implies that
\[ [s^{-n}\otimes l,s'^{-m}\otimes l']= s^{-n}d_{\tau(l)}(s'^{-m})\cdot l' + s'^{-m} [s^{-n}\otimes l, l'] \]
\[= s^{-n}d_{\tau(l)}(s'^{-m})\cdot l' - s'^{-m}d_{\tau(l')}(s^{-n})\cdot l + s'^{-m}s^{-n}[l,l']\]
Since this final expression only depends on the Lie bracket in $L$, and the extended anchor map, the Lie bracket on $L'$ is completely determined.
\end{proof}

The above technique for localizing Lie algebroids is clearly compatible with compositions of localizations, and defines a sheaf of Lie algebroids on $X$, for $X$ affine.  In the case of $X$ not affine, this local data may be sheafified; we will call any sheaf of Lie algebroids obtained this way a \textbf{Lie algebroid on $X$}.

For $X$ affine, and $X'$ an affine open subscheme, $\mathcal{U}_{X'}L'=\mathcal{O}_{X'}\otimes_X \mathcal{U}_XL=\mathcal{U}_XL\otimes_X \mathcal{O}_{X'}$.  This means that localizing enveloping algebras is the same on the left and on the right; so from now on we can refer to localizing them without refering to a side.  An $\mathcal{O}_X$-bimodule which has the property that left localization is isomorphic to right localization will be called \textbf{nearly central}; since it means that as a sheaf on $X\times X$, it is supported scheme-theoretically on the diagonal.

The universal enveloping algebra a non-affine Lie algebroid $(X,L)$ will be defined as the sheaf of algebras $\mathcal{U}_XL$ which is affine-locally the enveloping algebra of $(X,L)$.  Since enveloping algebras are nearly central, this is a quasi-coherent sheaf as both a left and right $\mathcal{O}_X$-module.

It is worth noting that, while the global sections of a Lie algebroid $(X,L)$ is again a Lie algebroid $(\Gamma(X),\Gamma(L))$, the global sections of $\mathcal{U}_XL$ is not necessarily the enveloping algebra of $(\Gamma(X),\Gamma(L))$.  For example, take the tangent bundle on $\mathbb{P}^1$.  The global Lie algebroid is $(\mathbb{C},\mathfrak{sl}_2)$ with trivial anchor map, but the global sections of $\mathcal{D}_{\mathbb{P}^1}$ is the algebra $\mathcal{U}\mathfrak{sl}_2/c$, where $c$ is the casimir element.

\subsection{Koszul Complexes.}

The sheaf of graded algebras $\widetilde{\mathcal{U}}$ can then be defined in the natural way, and is again a sheaf of nearly central algebras.  Since $\mathcal{U}^i$ is a nearly central $\mathcal{O}_X$-bimodule, so is $^*(\mathcal{U}^i)=J^i$.  The tensor algebra $T_XJ^1$ is then also nearly central, and so is $T_XJ^1/\langle J^2\rangle=\widetilde{\mathcal{U}}^{\perp}$.  Thus, $\widetilde{\mathcal{U}}^{\perp}$ also defines a sheaf of nearly central algebras on $X$.  From this, it follows that the left Koszul complex is also compatible with localization.
\begin{lemma}\label{lemma:koszul local}
Let $X'$ be an open affine subscheme of affine $X$, and $L'$ the localization of $L$.  Then $K^\bullet_{(X',L')}$ is the localization of $K^\bullet_{(X,L)}$.
\end{lemma}
\begin{proof}
On the level of terms of the complex, 
\begin{eqnarray*}
\mathcal{O}_{X'}\otimes_X \widetilde{\mathcal{U}_XL}(-i)\otimes_X (\widetilde{\mathcal{U}_XL}^{\perp i})^*
&=&\widetilde{\mathcal{U}_{X'}L'}(-i)\otimes_X (\widetilde{\mathcal{U}_XL}^{\perp i})^*\\
&=&\widetilde{\mathcal{U}_{X'}L'}(-i)\otimes_{X'} \mathcal{O}_{X'}\otimes_X (\widetilde{\mathcal{U}_XL}^{\perp i})^*\\
&=&\widetilde{\mathcal{U}_{X'}L'}(-i)\otimes_{X'} (\widetilde{\mathcal{U}_{X'}L'}^{\perp i})^*
\end{eqnarray*}
Note that the key is that the enveloping algebra is nearly central, and so localizing on the left localizes on the right.  Finally, it is immediate to show that the Koszul boundary is the correct one, because the Koszul boundary was defined in terms of multiplication in $\mathcal{U}_XL$, and localization is an algebra homomorphism.
\end{proof}

From this, and analogous observations, it can be shown that all the complexes and bicomplexes defined in Section \ref{section:koszul} are nearly central and compatible with localization.

\subsection{Projective Geometry.}

Define in the obvious way the module categories $mod(\mathcal{U}_XL)$ and $gr(\widetilde{\mathcal{U}_XL})$, which are $\mathcal{O}_X$-quasicoherent sheaves of left modules of the appropriate sheaf of algebras.  In either category, the $Hom$ set is naturally a $\mathcal{O}_X$-bimodule, and is a nearly central bimodule.  Therefore, localization is independent of side, and we can canonically define the sheafy $\mathcal{H}om$ as the sheaf of $\mathcal{O}_X$-bimodules on $X$ which is locally the corresponding $Hom$.

The quotient category $qgr(\widetilde{\mathcal{U}_XL})$ can be defined identically to the affine case.  The sheafy $\mathcal{H}om$ can also be defined in this case, by
\[ \mathcal{H}om_{qgr(\widetilde{\mathcal{U}})}(\pi M,\pi N):=\mathcal{H}om_{gr(\widetilde{\mathcal{U}})}(M,\omega\pi N)\]
Let $T$ denote $\oplus_{i=0}^n\pi\widetilde{\mathcal{U}}(-i)$.  Then $\mathcal{H}om_{qgr(\widetilde{\mathcal{U}})}(T,T)$ is defined locally by Theorem \ref{thm:end T}, and so it is the sheaf of nearly central algebras
\[\mathcal{H}om_{qgr(\widetilde{\mathcal{U}})}(T,T)=E^{op}
=\left(\begin{array}{ccccc}
\mathcal{O}_X & \mathcal{U}^1 & \mathcal{U}^2 & \cdots & \mathcal{U}^n \\
0 & \mathcal{O}_X & \mathcal{U}^1 & \cdots & \mathcal{U}^{n-1} \\
0 & 0 & \mathcal{O}_X & \cdots & \mathcal{U}^{n-2} \\
\vdots & \vdots & \vdots &\ddots & \vdots \\
0 & 0 & 0 & \cdots & \mathcal{O}_X \\
\end{array}\right)\]
and $\mathcal{E}xt^i_{qgr}(T,T)=0$ for $i>0$.  

Let $mod(E)$ denote the category of quasi-coherent sheaves of left $E$-modules on $X$.  Then $\mathbf{R}\mathcal{H}om_{qgr}(T,-)$ defines a functor from $D^b(qgr(\widetilde{\mathcal{U}}))$ to $D^b(mod(E))$.  We then have the non-affine version of the main theorem.
\begin{thm}
The functor $\mathbf{R}\mathcal{H}om_{qgr}(T,-)$ defines a quasi-equivalence of dg categories, and thus an equivalence of triangulated categories
\[D^b(qgr(\widetilde{\mathcal{U}}))\simeq D^b(mod(E))\]
\end{thm}
\begin{proof}
Since this functor is a quasi-equivalence on affine local subsets, the theorem will follow from effective descent for dg categories.  Specifically, the sheafy categories $D^b(qgr(\widetilde{\mathcal{U}}))$ and $D^b(mod(E))$ can be obtained as a homotopy limit of the affine local categories.  Since the functor is a locally a quasi-equivalence, it must be one in the limit.  See section 7.4 in \cite{BeilinsonDrinfeldQuantization}.
\end{proof}

\subsection{Example: Differential Operators on Flag Manifold.}
The Beilinson-Bernstein theorem is a powerful theorem which characterizes the category of $\mathcal{D}$-modules on a flag manifold.  Combining this with the Beilinson equivalence gives an equivalence of derived categories between the global differential operators on a flag manifold and a $\mathcal{O}_X$-coherent algebra.

Let $G$ be a semisimple Lie group over $\mathbb{C}$, and let $B$ be a Borel subgroup.  Let $G/B$ be the flag manifold, and let $\mathcal{D}$ denote the sheaf of differential operators on $G/B$. Then the Beilinson-Bernstein theorem says that
\begin{itemize}
\item The global sections of $\mathcal{D}$ are isomorphic to $\mathcal{U}\mathfrak{g}/Z$, where $Z$ is the center of $\mathcal{U}\mathfrak{g}$.
\item $G/B$ is $\mathcal{D}$-affine; that is, the global section functor \[\Gamma:mod(\mathcal{D})\rightarrow mod(\Gamma(\mathcal{D}))=mod(\mathcal{U}\mathfrak{g}/Z)\] is an equivalence.
\end{itemize}
The isomorphism is filtration-preserving, and so the above equivalence extends
\[gr(\widetilde{\mathcal{D}})\simeq gr(\widetilde{\mathcal{U}\mathfrak{g}/Z}),\;\;\;
qgr(\widetilde{\mathcal{D}})\simeq qgr(\widetilde{\mathcal{U}\mathfrak{g}/Z})\]
Combining this with the Beilinson equivalence, we get that
\[D^b(qgr(\widetilde{\mathcal{U}\mathfrak{g}/Z})) \simeq D^b(mod(E))\]
where $E$ is the sheaf of nearly central algebras
\[\left(\begin{array}{ccccc}
\mathcal{O}_X & \mathcal{D}^1 & \mathcal{D}^2 & \cdots & \mathcal{D}^d \\
0 & \mathcal{O}_X & \mathcal{D}^1 & \cdots & \mathcal{D}^{d-1} \\
0 & 0 & \mathcal{O}_X & \cdots & \mathcal{D}^{d-2} \\
\vdots & \vdots & \vdots & \ddots & \vdots \\
0 & 0 & 0 & \cdots & \mathcal{O}_X \\
\end{array}\right)\]
This provides yet another example of a graded algebra which is derived equivalent to one of these upper triangular algebras.



\bibliography{MyUniBib}{}
\bibliographystyle{plain}
%
%
%
%
%
%
%
%
%
%
%
%
%
%
%
%
%

\end{document}